# Quasi-Steady-State Approach for Efficient Multiscale Simulation and Optimization of mAb Glycosylation in CHO Cell Culture

Yingjie Ma, Jing Guo, Andrew J. Maloney[1] and Richard D. Braatz*

Massachusetts Institute of Technology, Cambridge, MA, 02139, USA

## Abstract

Glycosylation is a critical quality attribute for monoclonal antibody (mAb) production, influenced by both process conditions and cellular mechanisms. Multiscale mechanistic models, spanning from the bioreactor to the Golgi apparatus, have been proposed for analyzing the glycosylation process. However, these models are computationally intensive to solve when using traditional methods, making optimization and control challenging. In this work, we propose a quasi-steady-state (QSS) approach for efficiently solving the multiscale glycosylation model. By introducing the QSS assumption and assuming negligible nucleotide sugar donor (NSD) flux for glycosylation in the Golgi, the large-scale partial differential algebraic equation system is converted into a series of independent differential algebraic equation systems. Based on that representation, we develop a three-step QSS simulation method and further reduce computational time through parallel computing and nonuniform time grid strategies. Case studies in simulation, parameter estimation, and dynamic optimization demonstrate that the QSS approach can be more than 300-fold faster than the method of lines, with less than 1.6% relative errors. This work establishes a solid foundation for multiscale model-based optimization and control of the glycosylation process, supporting the implementation of quality by design.

---

* To whom correspondence should be addressed. braatz@mit.edu.

[1] Present address: Amgen, Cambridge, Massachusetts, USA.

## 1 Introduction

Monoclonal antibody (mAb) is the dominant product in the biopharmaceutical market in terms of both the number of approvals (53.5% biopharmaceutical approvals were mAbs from 2018 to 2022) and the commercial value (15 out of the top 20 biopharmaceutics by sales were mAbs in 2021) [1]. Because of their excellent specificity and affinity for both secreted and cell-surface targets, mAbs can be used to treat cancers, autoimmune diseases, inflammatory diseases and so on [2]. Although the biological activity of mAbs is encoded in the sequence of amino acids, it is also significantly influenced by post-translational modifications [3]. One of the most important post-translational modifications for mAbs is N-glycosylation, which can significantly impact the stability, immunogenicity, and efficacy of protein therapeutics [4–6]. Therefore, it is a critical quality attribute (CQA) for mAb therapeutics [3].

N-glycosylation refers to the cotranslational covalent attachment of an oligosaccharide group to an asparagine (Asn) side chain in secreted and membrane glycoproteins [4]. The process begins in the endoplasmic reticulum (ER) with the addition of a precursor nine-mannose oligosaccharide on the protein. After proper folding in the ER, the glycoprotein is transported to the Golgi apparatus, where a network of thousands of enzyme-catalyzed reactions modifies the precursor oligosaccharide, generating a large number of diverse glycan structures [7]. Therefore, as the products of a complex reaction network, the types and quantities of glycans are not readily revealed by the genetic code [8]. Instead, they are influenced by a multitude of reaction-related factors, such as the accessibility, concentrations, and kinetics of related enzymes in the cell, metabolite concentrations, and cell culture operating conditions (e.g., temperature, pH) [8–10]. Due to the complex interactions among these elements, tuning the glycan distribution through extensive experiments is both costly and time-consuming [10]. To address such issue, mathematical models were proposed to predict the glycosylation profile so that the glycan distribution can be optimized and controlled with less experiments [9,11–14]. In the remainder of this introduction, we first survey existing glycosylation models and then examine the state-of-the-art multiscale model. Although this model is accurate, its computational expense currently prevents its use in optimal control. To remove this barrier, the present work develops a strategy that sharply accelerates the simulation, making the model suitable for optimization-based control.

Data-driven glycosylation models have been proposed in the literature, such as response surface [15], Markov chain [8], and artificial neural network models [16]. However, considering the need for interpretability and extrapolation, it is desirable to develop mechanistic models for prediction of the glycan distribution [17]. Shelikoff et al. [11] proposed a model to predict the precursor oligosaccharide in the ER, while the other researchers mainly focus on developing models for the complicated glycosylation reactions in the Golgi apparatus. Umaña and Bailey [18] and Krambeck and Betenbaugh [12] modelled the Golgi apparatus as four serial continuous stirred tank reactors (CSTRs) to predict the glycans attached on mAbs. Hossler et al. [19] compared the results of modeling Golgi apparatus as four serial CSTRs and four serial plug flow reactors (PFRs), which showed that the later predicted glycan distribution better, supporting the Golgi maturation model. Also based on the cisternal maturation assumption, Jimenez del Val. et al. [13] approximated the Golgi apparatus by a single dynamic PFR model with the consideration of enzyme distribution along the Golgi length and the transport of NSDs into the Golgi apparatus. Their model presented improved predictive ability compared to previous models [12,19]. However, all the above models only modelled the ER or Golgi apparatus alone without its link with process conditions. To address that issue, Jedrzejewski et al. [20] proposed a multiscale model incorporating the cell culture submodel at the bioreactor level, nucleotide sugar donor (NSD) synthesis submodel at the intracellular level, and the glycosylation PFR model inside the Golgi. This multiscale model was later used for the optimization of the glycosylation process, generating an optimal operating strategy that improved galactosylation for 90% compared to the control experiment [9]. Villiger et al. [21] modified Jimenez del Val et al.'s Golgi model to consider the effects of pH and manganese on glycosylation, and then connected the Golgi model with the cell culture model for process simulation. The multiscale mechanistic model was shown to have superior extrapolation ability than a response surface model [22]. Although the dynamic PFR Golgi model was found to outperform the CSTR Golgi model, the resulting partial differential algebraic equation (PDAE) system can be 150-fold slower to simulate than the latter when traditional solution methods such as the method of lines (MOL) [21] are used [17], making its use in glycosylation process analysis cumbersome and even more challenging for optimization and control.

To solve the dynamic multiscale model more efficiently, the quasi-steady-state (QSS) approach has been used by several authors to simplify the Golgi PDAE model as a differential algebraic equation (DAE) model. The simplification assumes that the Golgi submodel is approximately at steady state, as the dynamics within the Golgi submodel are much faster than those in the outer level submodels. Jimenez del Val et al. [23] estimated parameters in the Golgi glycosylation submodel by using the standalone steady-state Golgi submodel. However, this model needed to use intracellular glycan data, which are rarely available and usually approximated by using extracellular glycan data, introducing further errors for parameter estimation. Maloney [24] applied the QSS method for the efficient simulation of a perfusion reactor described by a multiscale mAb glycosylation model, and assumed that extracellular glycosylation profiles are the same as intracellular glycosylation profiles. Nevertheless, the assumption may not hold for a dynamic process far from the steady state, especially for the fed-batch process. Overall, neither QSS method linked process operating conditions with dynamic extracellular glycan profiles in general operating conditions. Although after applying the QSS assumption, all the submodels in the multiscale model involve only DAEs, they are defined in different domains. The cell culture and NSD synthesis submodels are defined in the temporal domain, while the innermost Golgi model is in the spatial domain. Therefore, it is infeasible to solve the multiscale model by lumping all the three DAE submodels together and treating them as a single large DAE model for numerical integration. Actually, such a direct combination recreates a PDAE model, making the computation intensive again.

In the current work, we develop a novel QSS simulation method that can predict dynamic extracellular glycan distribution based on the multiscale glycosylation model [9]. Besides the QSS assumption, a second assumption of negligible NSD flux consumed in Golgi glycosylation is introduced to decouple the outer-level NSD submodel from the innermost Golgi glycosylation submodel. With that, we first simulate the dynamic cell culture submodel and NSD submodel to generate the dynamic trajectories of state variables for subsequent use as inputs to the Golgi glycosylation model. Sequentially, with above inputs, we perform steady-state simulations of the Golgi glycosylation submodel at various time points to obtain the intracellular glycan trajectories. These intracellular glycan compositions are then used as time-variant parameters in the cell culture submodel to compute the extracellular glycan trajectories. To further enhance the computational speed, we propose applying parallel computing and nonuniform time grid strategies to the QSS simulation. The comparison between PDE and QSS simulations shows that the proposed simulation method can be more than 300-fold faster with relative errors of less than 2%. Finally, the QSS simulation method is used in parameter estimation and dynamic optimization, which further validates its efficiency and accuracy.

The next section introduces the adopted multiscale glycosylation model. Then Section 3 develops the QSS simulation method step by step with assumptions and accuracy validated. Section 4 proposes to speed up the QSS simulation by parallel computing and nonuniform time points allocation. Section 5 applies the QSS simulation to parameter estimation and dynamic optimization problems. Finally, conclusions are provided in Section 6.

## 2 Glycosylation model

The multiscale model is mainly based on the model from Kotidis et al. [9], while the effects of manganese and ammonia on glycosylation are also considered according to Villiger et al. [21] The overall model can be illustrated in Fig. 1, which includes cell culture, NSD, and Golgi submodels. In the following part of this section, we give a brief introduction to the model, while the detailed model equations are shown in Supplementary Material S1.

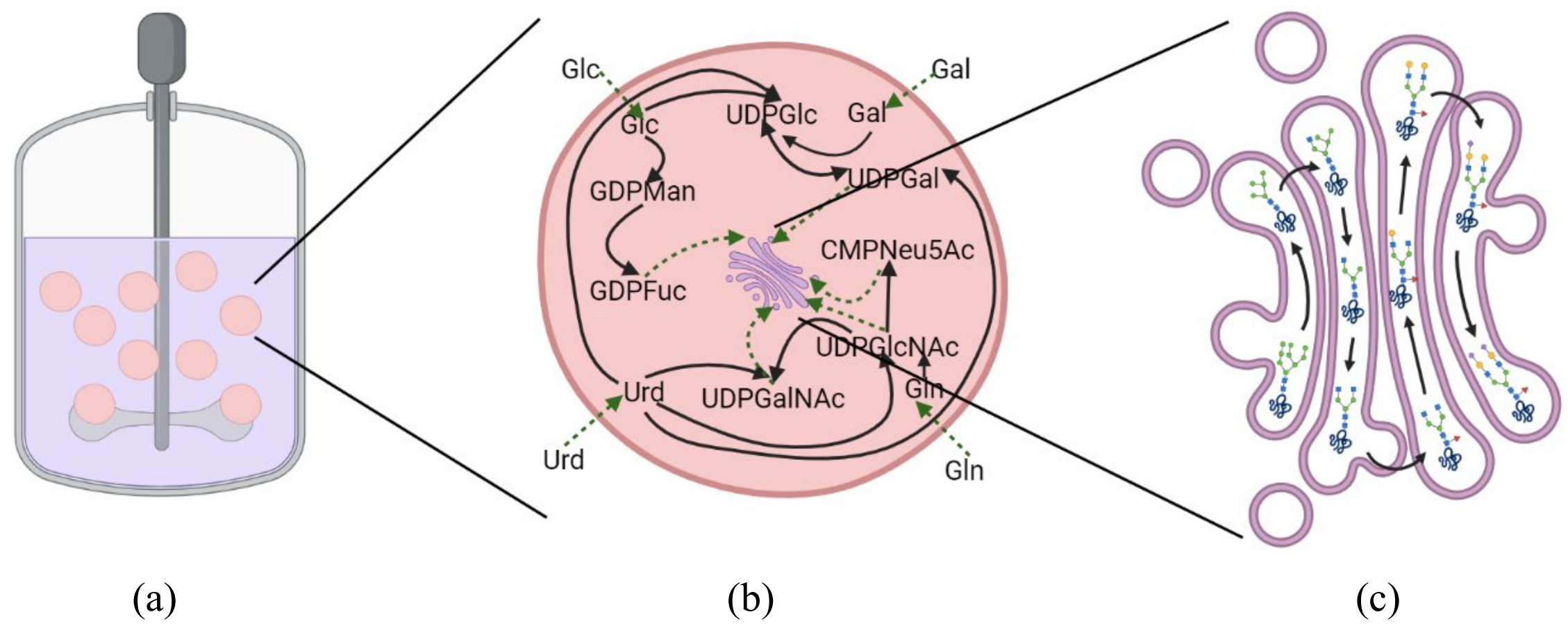

Figure 1. Multiscale glycosylation model involving (a) cell culture model, (b) NSD synthesis model, (c) Golgi glycosylation model.

At the bioreactor level, an unstructured differential algebraic equation (DAE) model is used to describe the cell growth, death, and metabolism. The metabolites considered in the model include ammonia (Amm), asparagine (Asn), aspartate (Asp), glucose (Glc), galactose (Gal), glutamine (Gln), glutamate (Glu), lactose (Lac), and uridine (Urd). Glucose and asparagine are regarded as the limiting substrates, while lactate, ammonia and uridine are inhibiting metabolites. Ammonia and uridine can cause the death of cells. The mass balance equations for cell numbers and metabolites are

$$\frac{d(VX_v)}{dt} = (\mu - \mu_{\text{death}})VX_v - F_{out}X_v, \tag{1}$$

$$\frac{d(V[\text{Metabolite}])}{dt} = F_{\text{in}}[\text{Metabolite}_{\text{in}}] - F_{\text{out}}[\text{Metabolite}] + q_{\text{metabolite}}VX_v, \tag{2}$$

where $V$ is the cell culture volume (L), $\mu$ ($\mathrm{h^{-1}}$) is the specific cell growth rate, $\mu_{\mathrm{death}}$ ($\mathrm{h^{-1}}$) is the specific cell death rate, $X_v$ ($\mathrm{cell \cdot L^{-1}}$) is the viable cell density, $F_{\mathrm{in}}$ ($\mathrm{L \cdot h^{-1}}$) is the inlet flow rate of the bioreactor, and $F_{\mathrm{out}}$ ($\mathrm{L \cdot h^{-1}}$) is the outlet flow rate of the bioreactor. $q_{\mathrm{metabolite}}$ is the reaction rate of each metabolite, which has the unit $\mathrm{pg \cdot cell^{-1} \cdot h^{-1}}$ for mAb and the unit $\mathrm{mmol \cdot cell^{-1} \cdot h^{-1}}$ for the other metabolites. Their computation is shown in Eqs. (S9–S19), while the metabolic network is shown in Fig. S1. $[\mathrm{Metabolite_{in}}]$ (mM) and $[\mathrm{Metabolite}]$ (mM) are the concentrations of each metabolite in the inlet stream and bioreactor, respectively, which has the unit $\mathrm{pg \cdot L^{-1}}$ for mAb and the unit mM for the other metabolites.

At the cellular level, the NSD submodel is a DAE system describing the synthesis of NSDs, which is influenced by the extracellular glucose, galactose, uridine, and glutamine. The simplified NSD metabolic pathway is shown in Fig. 1b. The NSD submodel consists of mass balance equations,

$$\frac{d\left[\mathrm{NSD}_i^{\mathrm{intra}}\right]}{dt} = \sum_{j=1}^{N_{R1}} v_{i,j}^{\mathrm{nsd}} r_j^{\mathrm{nsd}} - f_{\mathrm{NSD}_i}^{\mathrm{hcp/lipid}} - f_{\mathrm{NSD}_i}^{\mathrm{precursor}} - f_{\mathrm{NSD}_i}^{\mathrm{glyc}}, \; i = 1,2,\cdots, N_{\mathrm{NSD}}, \tag{3}$$

for seven NSDs (GDPMan, GDPFuc, UDPGalNAc, UDPGlcNAc, CMPNeu5Ac, UDPGal and UDPGlc) and a series of Michaelis-Menten saturation kinetics for reaction rate computation. Here, $[\mathrm{NSD}_i^{\mathrm{intra}}]$ (mM) is the concentration of $\mathrm{NSD}_i$ in the cytosol, $N_{R1}$ is the total number of reactions producing or consuming NSDs, $v_{i,j}^{\mathrm{nsd}}$ ($-$) is the stoichiometric coefficient of $\mathrm{NSD}_i$ in the reaction $j$, $r_j^{\mathrm{nsd}}$ ($\mathrm{mmol \cdot L^{-1} \cdot h^{-1}}$) is the rate of reaction $j$, $f_{\mathrm{NSD}_i}^{\mathrm{hcp/lipid}}$ ($\mathrm{mmol \cdot L^{-1} \cdot h^{-1}}$) is the flux of $\mathrm{NSD}_i$ used for the synthesis of host cell proteins (hcp) and glycolipids, while $f_{\mathrm{NSD}_i}^{\mathrm{precursor}}$ ($\mathrm{mmol \cdot L^{-1} \cdot h^{-1}}$) and $f_{\mathrm{NSD}_i}^{\mathrm{glyc}}$ ($\mathrm{mmol \cdot L^{-1} \cdot h^{-1}}$) are the $\mathrm{NSD}_i$ demands for the precursor oligosaccharide formation and N-linked glycosylation, respectively. The latter is given by

$$f_{\mathrm{NSD}_i}^{\mathrm{glyc}} = \frac{[\mathrm{NSD}_i]}{K_{\mathrm{TP}_{\mathrm{NSD}_i}} + [\mathrm{NSD}_i]} r_{\mathrm{NSD}_i}^{\mathrm{glyc}}, \tag{4}$$

where $r_{\mathrm{NSD}_i}^{\mathrm{glyc}}$ ($\mathrm{mmol} \cdot \mathrm{L}^{-1} \cdot \mathrm{h}^{-1}$) is the consumption rate of $\mathrm{NSD}_i$ used for N-linked glycosylation in Golgi, and $K_{\mathrm{TP}_{\mathrm{NSD}_i}}$ (mM) is the transport protein saturation constant. $r_{\mathrm{NSD}_i}^{\mathrm{glyc}}$ is obtained from the following Golgi model, while the formula calculating $f_{\mathrm{NSD}_i}^{\mathrm{hcp/lipid}}$ and $f_{\mathrm{NSD}_i}^{\mathrm{precursor}}$ can be found in Supplementary Material S1.2.

NSDs supply sugars for the glycosylation reactions in the Golgi, and their concentrations in Golgi are determined by their concentrations in the cytosol. The dynamic Golgi model consists of reaction rate equations and mass balance equations of oligosaccharides. The glycosylation reaction network is shown in Fig. S2, while the expressions for reaction rates are given in Eqs. (S39–S44), while the mass balance equations are partial differential equations (PDEs) given by

$$\frac{\partial[\mathrm{OS}_i]}{\partial t} = -\mathrm{Vel}_{\mathrm{golgi}} \frac{\partial[\mathrm{OS}_i]}{\partial z} + \sum_{j=1}^{N_{R2}} \nu_{i,j} r_j , \qquad i = 1,2,\cdots,N_{\mathrm{OS}} \tag{5}$$

where $[\mathrm{OS}_i]$ is the concentration (μM) of oligosaccharide $i$ ($\mathrm{OS}_i$); $\mathrm{Vel}_{\mathrm{golgi}}$ ($\mathrm{Golgi\ length} \cdot \mathrm{min}^{-1}$) is the normalized linear velocity of the molecule passing through the Golgi, which is obtained from $q_{\mathrm{mab}}$ and the Golgi volume; $r_j$ ($\mathrm{\mu mol} \cdot \mathrm{L}^{-1} \cdot \mathrm{min}^{-1}$) is the reaction rate for reaction $j$; and $\nu_{i,j}$ is the stoichiometric coefficient of $\mathrm{OS}_i$ in reaction $j$. There are $N_{R2}$ reactions in total, and $N_{\mathrm{OS}}$ is the total number of oligosaccharides.

The reaction network in the current model is from Villiger et al. [21], which involves 33 oligosaccharides, 43 reactions, and 7 enzymes. Depending on the enzymes, there are three kinetics types for the reaction rate prediction, including Michaelis-Menten kinetics, sequential order Bi-Bi kinetics, and random order Bi-Bi kinetics as shown in the Supplementary Material. According to literature reports [25,26], manganese and ammonia have significant influence on glycosylation, so their effects are also included in our glycosylation kinetic equations as with Villiger et al. [21].

The Golgi glycosylation model predicts the glycan concentrations inside the Golgi apparatus. However, for real-world applications, we often need the extracellular glycoprotein concentrations in the cell culture media. The mass balance equations

$$\frac{d(V[\mathrm{GLY}_i^{\mathrm{extra}}])}{dt} = -F_{\mathrm{out}}[\mathrm{GLY}_i^{\mathrm{extra}}] + Vq_{\mathrm{mab}}[X_v]Y_i^{\mathrm{intra}}, \tag{6}$$

$$Y_i^{\mathrm{intra}} = \frac{[\mathrm{GLY}_i^{\mathrm{intra}}]}{[\mathrm{mAb}]}, \tag{7}$$

are used to compute the extracellular glycoprotein concentrations, where $Y_i^{\mathrm{intra}}$ $(\mathrm{mM} \cdot \mathrm{mM}^{-1})$ and $[\mathrm{GLY}_i^{\mathrm{intra}}]$ (mM) are the percentage and concentration of $\mathrm{GLY}_i$ glycosylated mAb that leaves the Golgi apparatus and enters the cytosol. There are six glycans considered – Man5, FA2G0, FA2G1, FA2G2, G0, and G2 – which consist of different sets of oligosaccharides. $[\mathrm{GLY}_i^{\mathrm{extra}}]$ (mM) is the extracellular concentration of mAb with $\mathrm{GLY}_i$ attached. The percentage of extracellular mAb attached with $\mathrm{GLY}_i$, $Y_i^{\mathrm{extra}}$ $(\mathrm{mM} \cdot \mathrm{mM}^{-1})$, the percentage of extracellular mAb attached with $\mathrm{GLY}_i$, is given by

$$Y_i^{\mathrm{extra}} = \frac{[\mathrm{GLY}_i^{\mathrm{extra}}]}{[\mathrm{mAb}]}. \tag{8}$$

Finally, the multiscale glycosylation model involving the above three submodels is a large-scale partial differential algebraic equation (PDAE) model with 30 ordinary differential equations (ODEs), 34 PDEs, and many strongly nonlinear algebraic equations (mainly kinetic equations). Although a dynamic simulation takes only a few seconds, sensitivity computation can take more than 10 minutes, depending on the number of dependent and independent variables, as shown later. This leads to an extended computation time that may exceed a day for dynamic optimization, making it impractical for nonlinear model predictive control and closed-loop model-based design of experiments. Therefore, we will develop the QSS method in the next section to significantly accelerate the simulation.

## 3. QSS Simulation for the Multiscale Glycosylation Model

This section states the assumptions for the QSS simulation, and validates the assumptions by numerical simulations under five different experimental operating strategies from Kotidis et al. [9]. All the experiments were conducted in shaking flasks with a working volume of 100 mL and run for 12 days or when cell viability was lower than 60%. In all the experiments, there were nutrient supplements every two days since day 2, and the feeding volume is 10% of the working volume. The nutrient compositions in the basic supplements are shown in Table S3, and the detailed experiment operating strategies are shown in Tables S4-S8. Manganese (II) chloride solution with a concentration of 1 μM was supplemented to the cell cultures at seeding. The difference of the five experiments lies in the concentrations of galactose and uridine in the nutrient supplements on certain days, which are shown in Table 1.

Table 1. Galactose and uridine addition in five different experiments.

| Experiment name | Galactose and uridine addition |
|---|---|
| Control | No galactose and uridine addition |
| 10G | Adding 10 mM galactose on days 4 and 8 |
| 10G5U | Adding 10 mM galactose and 5 mM uridine on days 4 and 8. |
| 10G20U | Adding 10 mM galactose and 20 mM uridine on days 4 and 8. |
| 50G5U | Adding 50 mM galactose and 5 mM uridine on days 4 and 8. |

In the paper, we will only show the results for the experiment 10G to save space, while the results for the other four experiments are shown in the Supplementary Material.

### 3.1 Time scales of different submodels and QSS assumption

First we make an assumption that substantially reduces the model complexity.

**Assumption 1** (QSS assumption): the Golgi glycosylation model is approximately at steady state for the given environment variables that can influence glycosylation reactions, which include extracellular ammonia concentration ($[\mathrm{Amm}]$), mAb flux entering Golgi ($q_{\mathrm{mab}}$), and intracellular NSD concentrations ($[\mathrm{NSD}_i^{\mathrm{intra}}]$).

This assumption reduces the PDEs (5) to the ODEs,

$$\frac{q}{V}\frac{\partial[\mathrm{OS}_i]}{\partial z} = \sum_{j=1}^{N_{R2}} \nu_{i,j} r_j \,, \qquad i = 1,2,\cdots,N_{\mathrm{OS}} \tag{9}$$

In the following, all the environment variables are denoted by the vector

$$\mathrm{env} = \left[[\mathrm{Amm}], q_{\mathrm{mab}}, \left[\mathrm{NSD}_i^{\mathrm{intra}}\right]\right]. \tag{10}$$

The QSS assumption is valid because the residence time ( $\tau = \frac{V}{q} < 30$ minutes) of oligosaccharides in the Golgi apparatus is much shorter than the time scale (several days) of the environment variables (not considering the event time points). Hence, env can be treated as constants within $\tau$ as shown in Fig. 2. Moreover, according to the property of the PFR model, the glycan distributions inside the Golgi are determined by env except for the initial $\tau$, and the distributions reach to steady state within $\tau$ for constant env. Because of these two factors, at each time point after the initial $\tau$, the glycan distributions in the Golgi model are approximately at steady state, i.e., quasi-steady state.

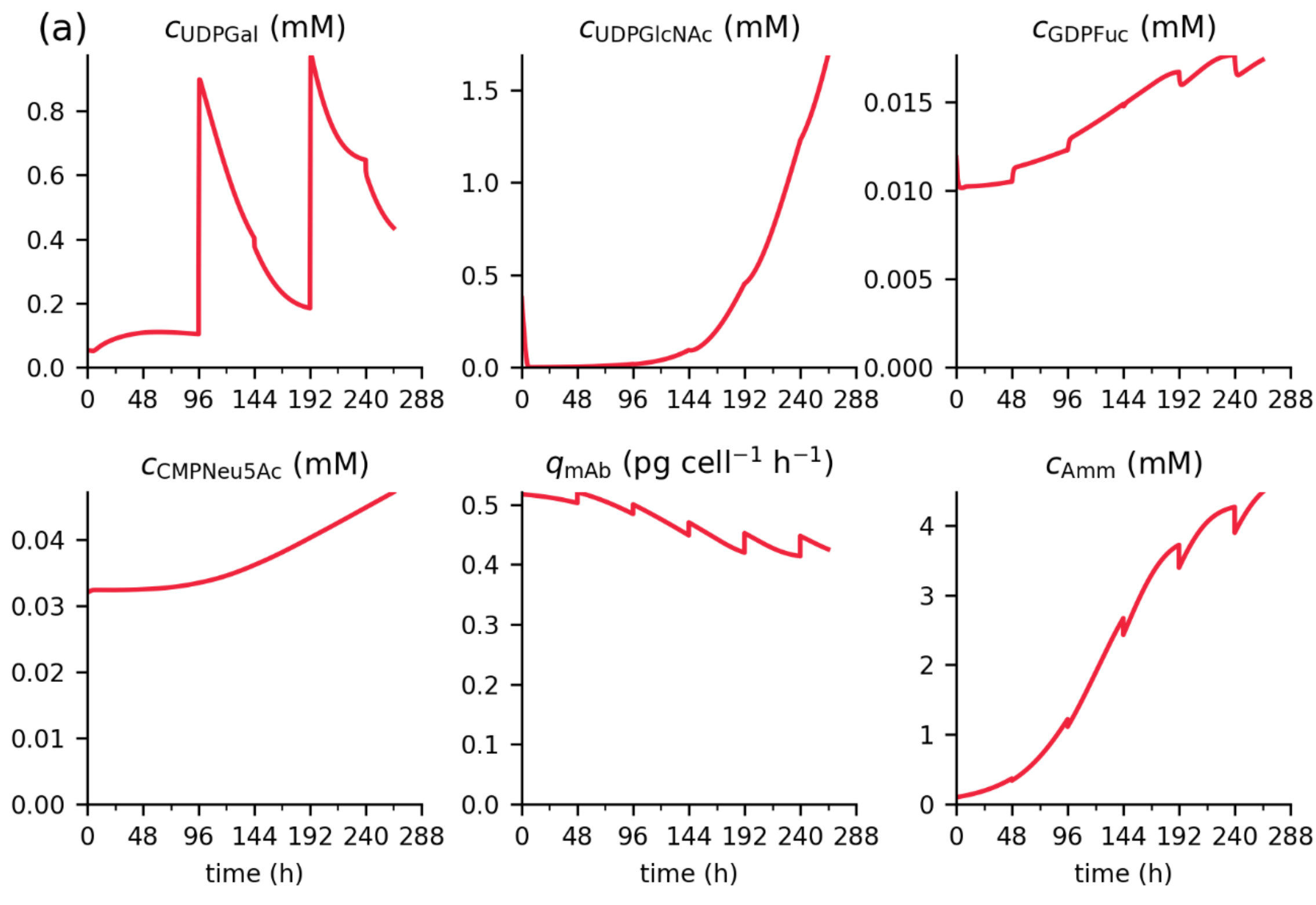

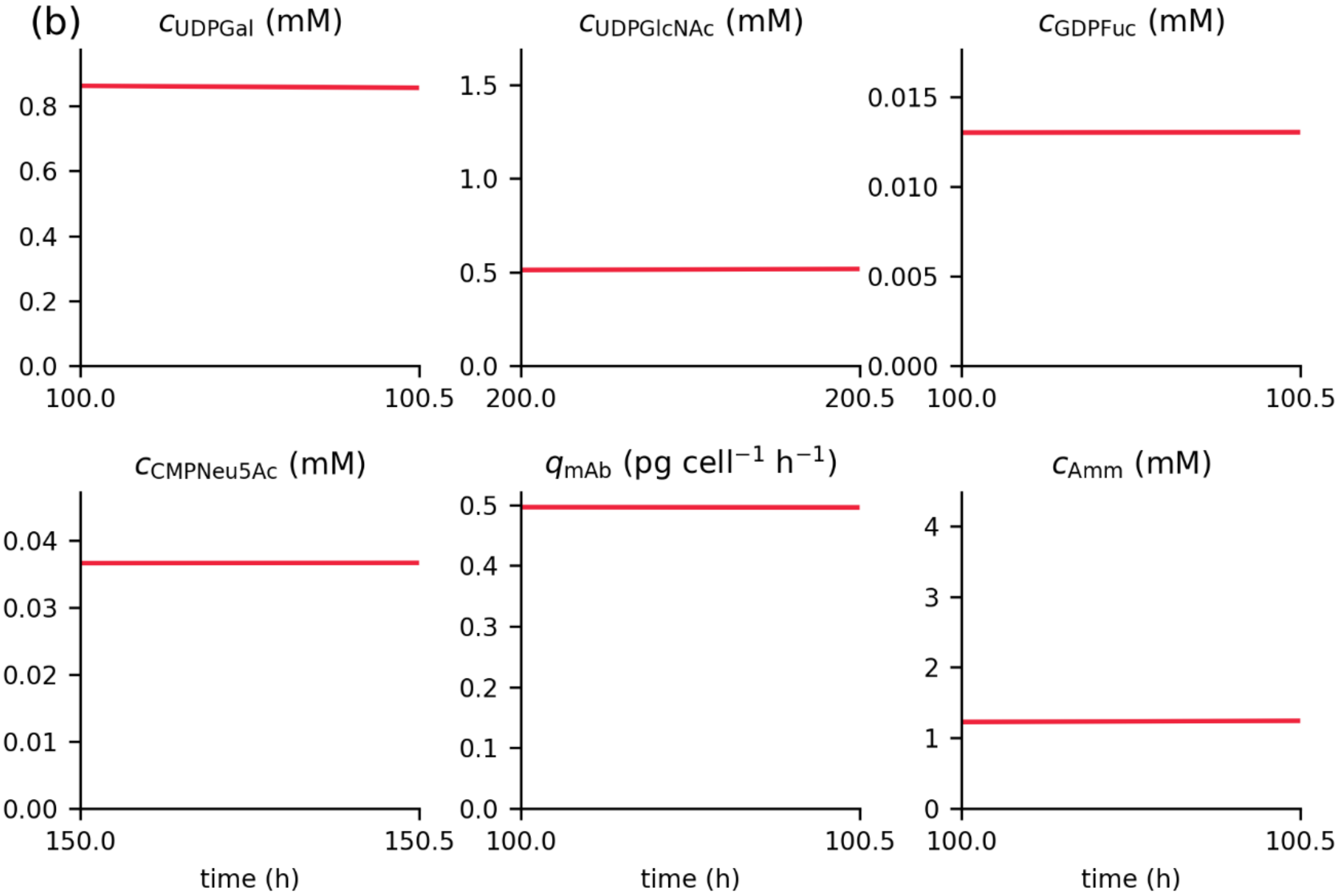

Figure 2. The (a) whole, and (b) half an-hour zoom-in trajectories of the environment variables. The variables include critical NSD concentrations in the cytosol, the flow rate of mAb entering the Golgi, and the ammonia concentration in the bioreactor.

**Remark 1**: The QSS assumption may not hold during the initial $\tau$ period of cell culture; however, this results in negligible error when predicting the accumulated extracellular glycan concentrations after a few hours in the QSS simulation. This is because $\tau$ is small, and the cell density during this period is very low, both contributing to negligible accumulation of glycoproteins in the bioreactor during the initial τ.

**Remark 2**: The short time scale (< 30 minutes) of the Golgi glycosylation model is due to the brief residence time of glycoproteins in the Golgi apparatus according to the property of the PFR model, rather than the glycosylation reaction rates. Since the residence time is consistently short, the QSS assumption remains valid for the multiscale glycosylation model.

The slow changes of the NSD concentrations result from the slow changes of the variables in the cell culture submodel, but the dynamics in the NSD submodel itself are actually fast. This can be seen from the jump of the concentrations of UDPGal, UDPGlcNAc, and GDPFuc at certain time points with feeding and/or at $t = 0$, where the NSD concentrations are far from the steady states at those time points. However, the dynamics of the NSD submodel are slower than those in the Golgi model, which has a time scale ranging from half an hour to a few hours (Fig. 3). As seen in Fig. 3, the concentrations of UDPGal and GDPFuc march to new quasi-steady states about 0.5 and 1 hour after the feeding of galactose, respectively. The concentration of UDPGlcNAc has dramatic change only at the beginning of the cell cultivation, and it takes around 6 hours to reach to its quasi-steady state. The influence of the times scales of the NSD submodel on QSS simulations are discussed in Sections 3.3 and 4.2.

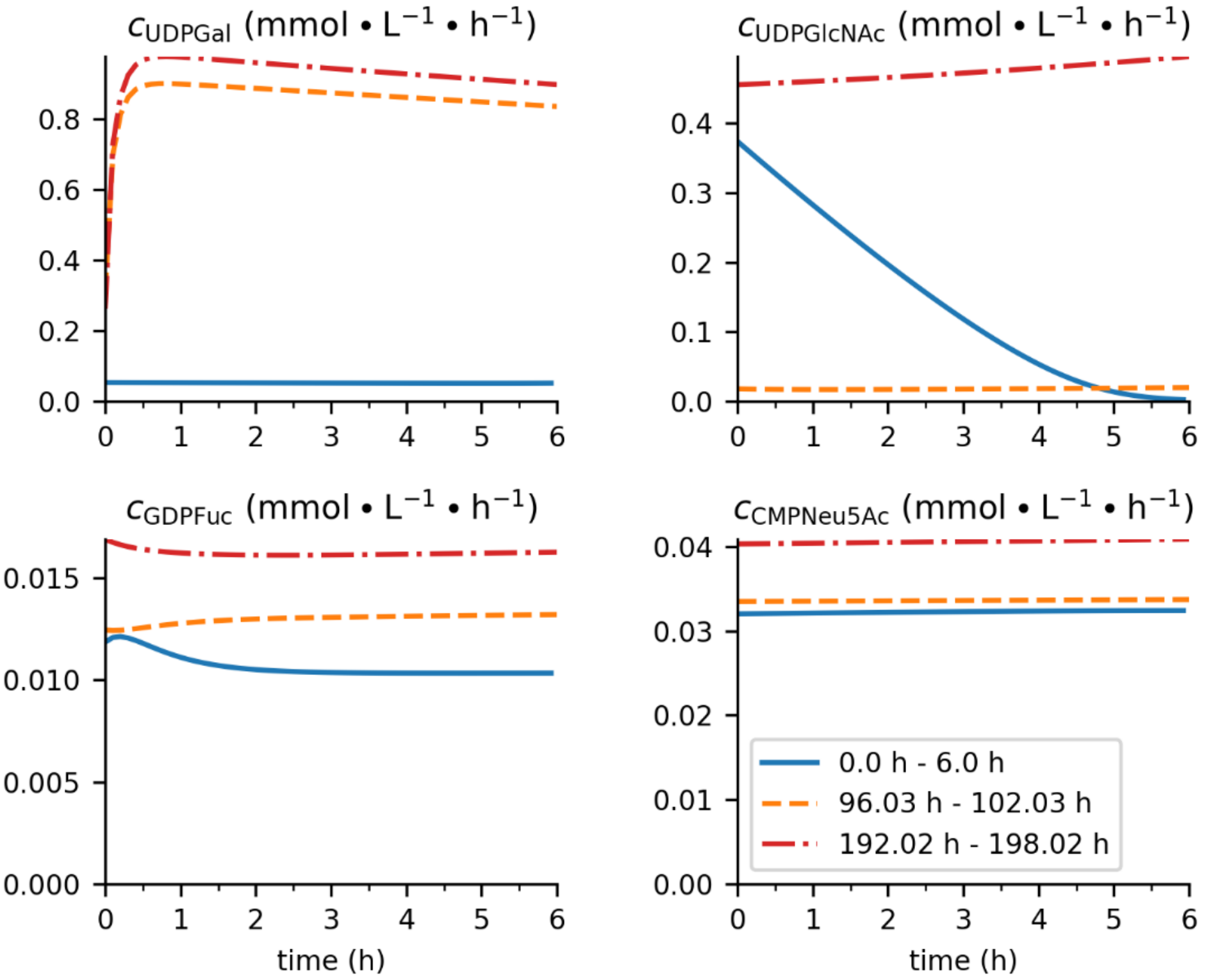

Figure 3 NSD concentration trajectories across three time intervals: after the initial cultivation period (0.00 h – 6.00 h), and following galactose feeding at 96.03 h – 102.03 h and 192.02 h – 198.02 h.

### 3.2 Negligible NSD flux consumed by glycosylation in the Golgi

Even after simplifying the dynamic Golgi model to a steady-state model according to the QSS assumption, we still cannot simulate the multiscale glycosylation model without solving PDEs. This is because the NSD submodel and the Golgi submodel are coupled and must be solved simultaneously, while they are defined in temporal domain and spatial domain, respectively. Therefore, it is infeasible to combine the three DAE submodels and solve them directly by one single DAE integration. To avoid the discretization in the spatial domain for the dynamic simulation of the multiscale glycosylation process, it is necessary to decompose the simulations in the temporal domain and spatial domain by decoupling the NSD model and the Golgi model. The coupling between the two submodels exists because the NSD model provides NSD concentrations for the Golgi model, while the NSD fluxes consumed in glycosylation, $f^{\text{glyc}}_{\text{NSD}_i}$ appearing in Eq. (3) of the former, require the computational results from the later. By using Assumption 2, the two submodels can be decoupled.

**Assumption 2**: The NSD fluxes consumed by glycosylation in the Golgi are negligible compared to the other NSD fluxes. That is, ignoring this term in the mass balance equation (3) for the intracellular NSD submodel has no significant impact on the prediction of intracellular NSD concentrations.

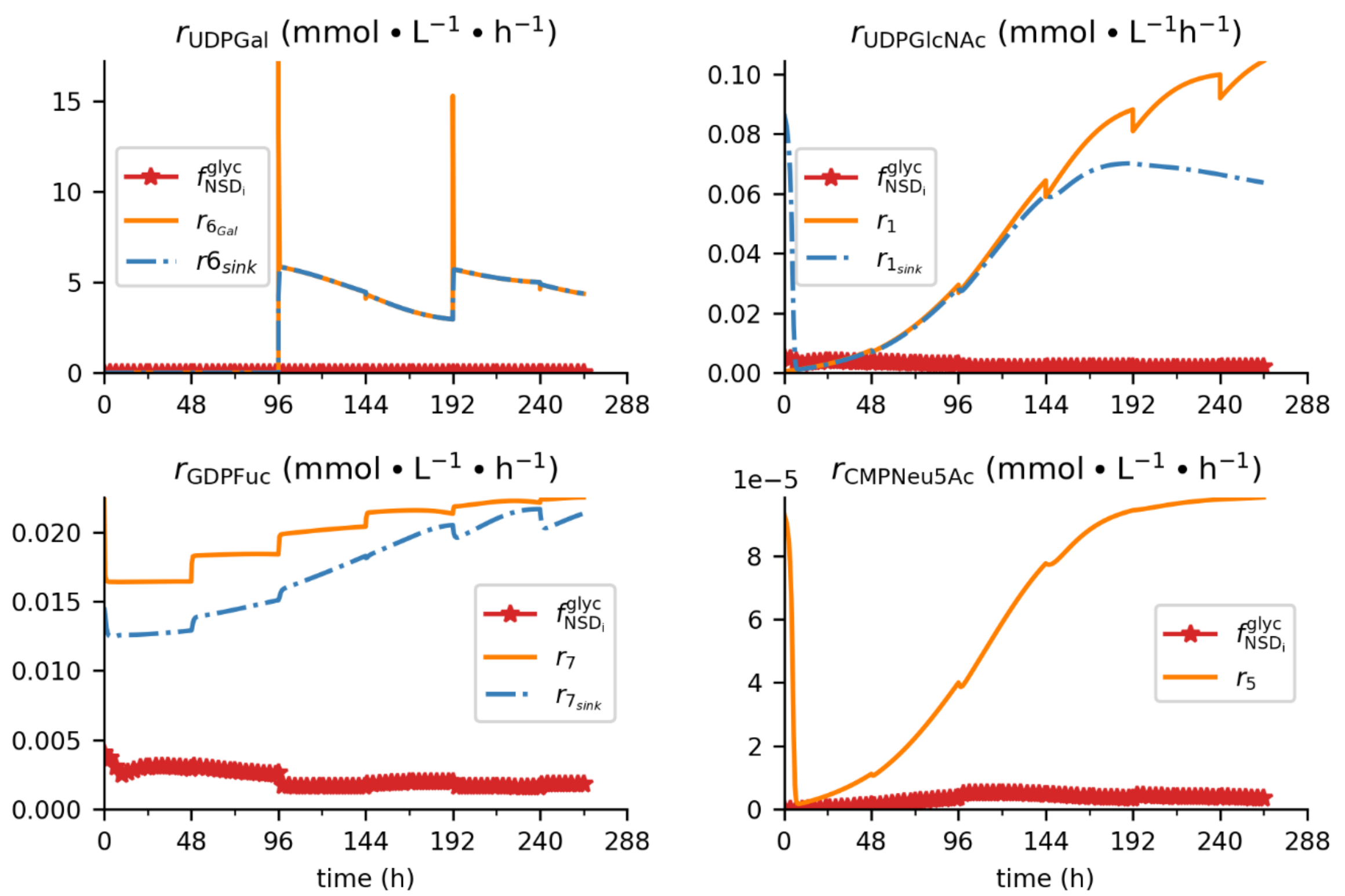

Figure 4. Comparison of NSD fluxes. nsd_flux_golgi refers to the NSD fluxes used for the glycosylation in the Golgi; the other fluxes starting with "r" refer to the reaction rates for the reactions in the cytosol.

The NSD fluxes for the glycosylation ($f^{\text{glyc}}_{\text{NSD}_i}$) in the Golgi and the other NSD fluxes are compared in Fig. 4, which indicates that the formers are orders of magnitude smaller than the reaction fluxes ($r_1$, $r_{1_{\text{sink}}}$, $r_7$, $r_{7_{\text{sink}}}$, $r_{6_{\text{Gal}}}$, $r_{6_{\text{sink}}}$, and $r_5$) in the cytosol. Therefore, the predicted NSD concentration trajectories show minimal changes after neglecting $f_{\text{NSD}_i}^{\text{glyc}}$ in Eq. (3), as demonstrated in Fig. 5, where two NSD concentration trajectories from the PDE simulations with and without considering $f_{\text{NSD}_i}^{\text{glyc}}$ nearly overlap.

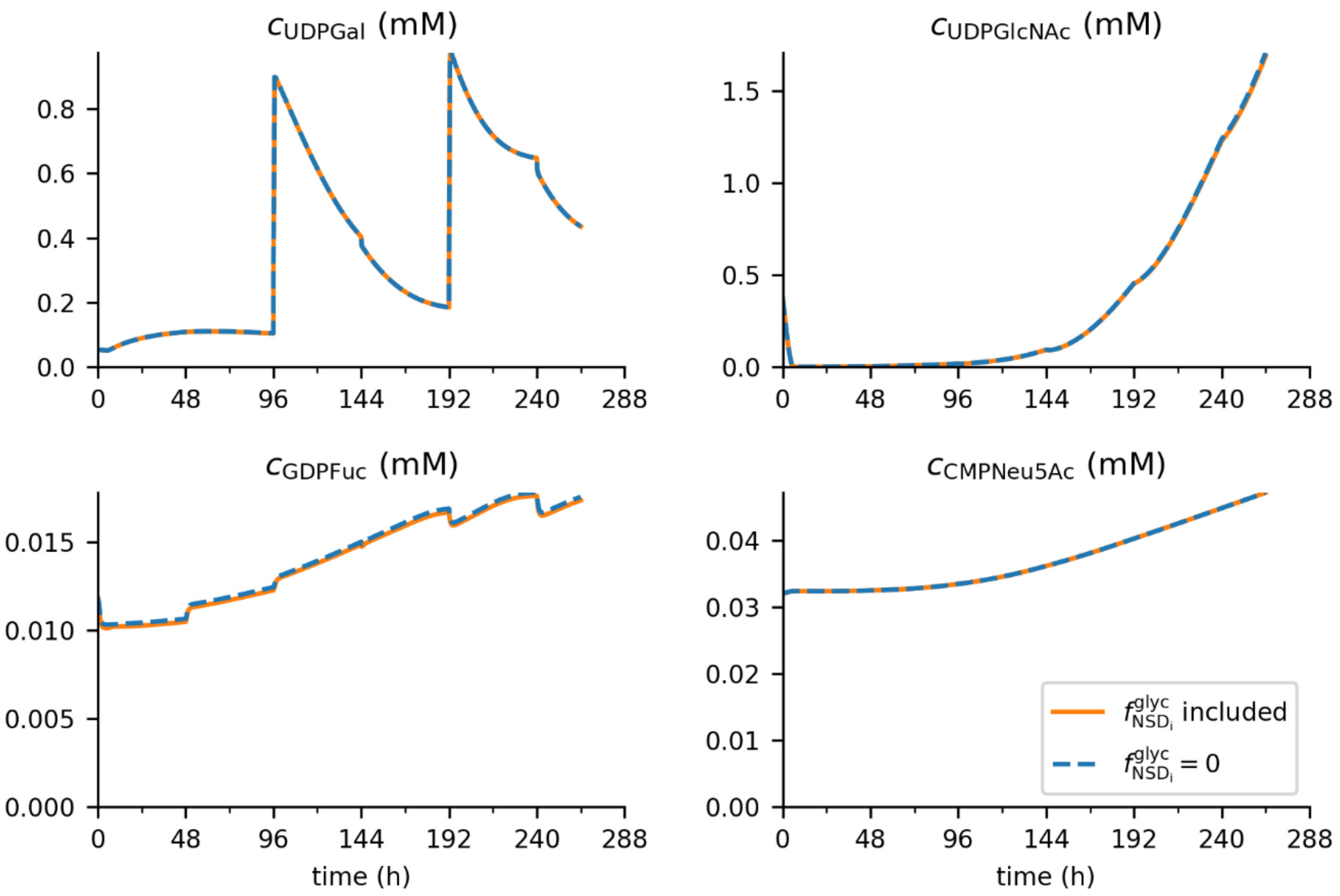

Figure 5. The trajectories of intracellular NSD concentrations from the PDE simulations with and without considering the NSD fluxes for the glycosylation in Golgi.

### 3.3 QSS Simulation for the intracellular glycan trajectory

Based on Assumption 2, we can neglect the term $f_{\text{NSD}_i}^{\text{glyc}}$ in the NSD model and simulate the cell culture model and NSD model first to generate env at a set of time points $\mathcal{T} \coloneqq \{t_k | k = 0,1,2,\cdots,K\}$. Then a series of steady-state simulations for the Golgi glycosylation submodel at time points $\mathcal{T}$ are conducted, and the concatenation of the simulation results at those time points approximate the trajectories of intracellular glycans, as shown in Fig. 6.

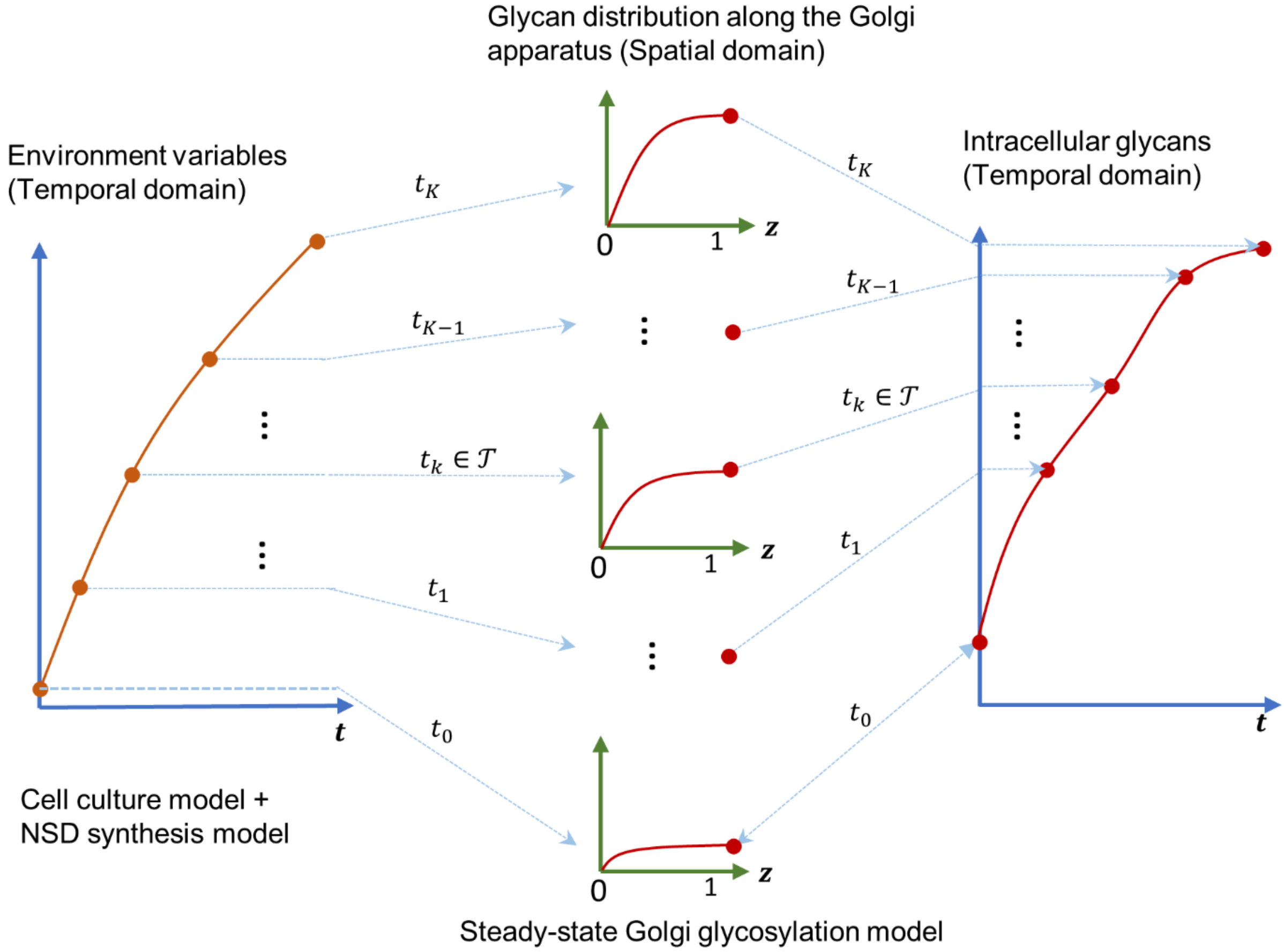

Figure 6. The computation of intracellular glycan trajectories by the QSS simulation approach.

When constructing the trajectories of intracellular glycan concentrations and percentages, there are three considerations for the selection of $t_k \in \mathcal{T}$:

(1) $\mathcal{T}$ should at least include all the time points right before and after the events where there might be significant change of $[\mathrm{GLY}_i^{\mathrm{intra}}]$ and $Y_i^{\mathrm{intra}}$ caused by the sudden change of env, e.g., the pulse feed of galactose and/or uridine. In the following part, we use $\mathcal{E}:=\{e|1,2,\cdots E\}$ to represent the set of such events, and the corresponding time set is $\mathcal{T}_e := \{t_e|e \in \mathcal{E}\}$. The time right after $t_e$ is notated as $t_e^+$, and we define the set which is $\mathcal{T}_e^+ := \{t_e^+|e \in \mathcal{E}\}$. Here, $t_e^+ = t_e + \epsilon$, and $\epsilon$ is a very small value, e.g. 0.01 h, so that both the computational time and the accumulated glycans within $[t_e, t_e^+]$ are negligible even if there might be a jump of $Y_i^{\mathrm{intra}}$ in the interval.

(2) Besides $\mathcal{T}_e$ and $\mathcal{T}_e^+$, the time points where $c_{\mathrm{NSD}_i}$ change dramatically may need to be included in $\mathcal{T}$ to capture the trajectory of $Y_i^{\mathrm{intra}}$, which is important for the computation of extracellular glycanform. These points depend on the time scales of the NSD submodel dynamics.

(3) The inclusion of more time points in $\mathcal{T}$ can potentially obtain better approximation to the trajectories of $[\mathrm{GLY}_i^{\mathrm{intra}}]$ and $Y_i^{\mathrm{intra}}$, but it will also increase the computational time because of more DAE simulations required for the steady-state Golgi model.

With the above considerations, the QSS simulations using two different sets of time points are used to demonstrate the method. The first set ($\mathcal{T}_{\mathrm{event}}$) includes only event-related time points $t_e$ and $t_e^+$, while the second set ($\mathcal{T}_{100}$) consists of all the above time points and additional 100 uniformly distributed time points in the cell culture period to satisfy the requirements of (2) and (3), i.e.,

$$\mathcal{T}_{\mathrm{event}} = \mathcal{T}_e \cup \mathcal{T}_e^+, \tag{11}$$

$$\mathcal{T}_{100} = \mathcal{T}_{\mathrm{event}} \cup \left\{ t_k | t_k = k \frac{T}{100-1}, k = 1, 2, \cdots, 100 \right\}, \tag{12}$$

where $T$ is the cell culture period. We have also tried to add 200 additional time points for the QSS simulation, but there is no noticeable change in the simulation results, so those simulations will not be discussed in the paper. In the simulation for experiment 10G, $\mathcal{T}_{\mathrm{event}}$ has 28 time points, while $\mathcal{T}_{100}$ has 116 time points, which can be visualized as Fig. 7.

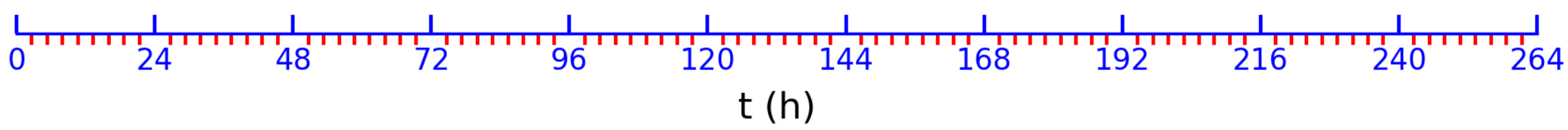

Figure 7. The time points in $\mathcal{T}_{\mathrm{event}}$ and $\mathcal{T}_{100}$. The long ticks denote the time points in $\mathcal{T}_{\mathrm{event}}$, while the short ticks denote the time points newly added in $\mathcal{T}_{100}$. The time points in $\mathcal{T}_e^+$ and $\mathcal{T}_e$ are so close that they are indistinguishable in the figure.

The QSS simulation using $\mathcal{T}_{\mathrm{event}}$ is notated as QSSEvent, while the QSS simulation using $\mathcal{T}_{100}$ is notated as QSS100. The comparison of $Y_i^{\mathrm{intra}}$ from the two QSS simulations and a PDE simulation (notated as PDE400) are shown in Fig. 8. Here, the PDE simulation uses the MOL method and apply the finite difference method in the spatial domain using 400 discretization points.

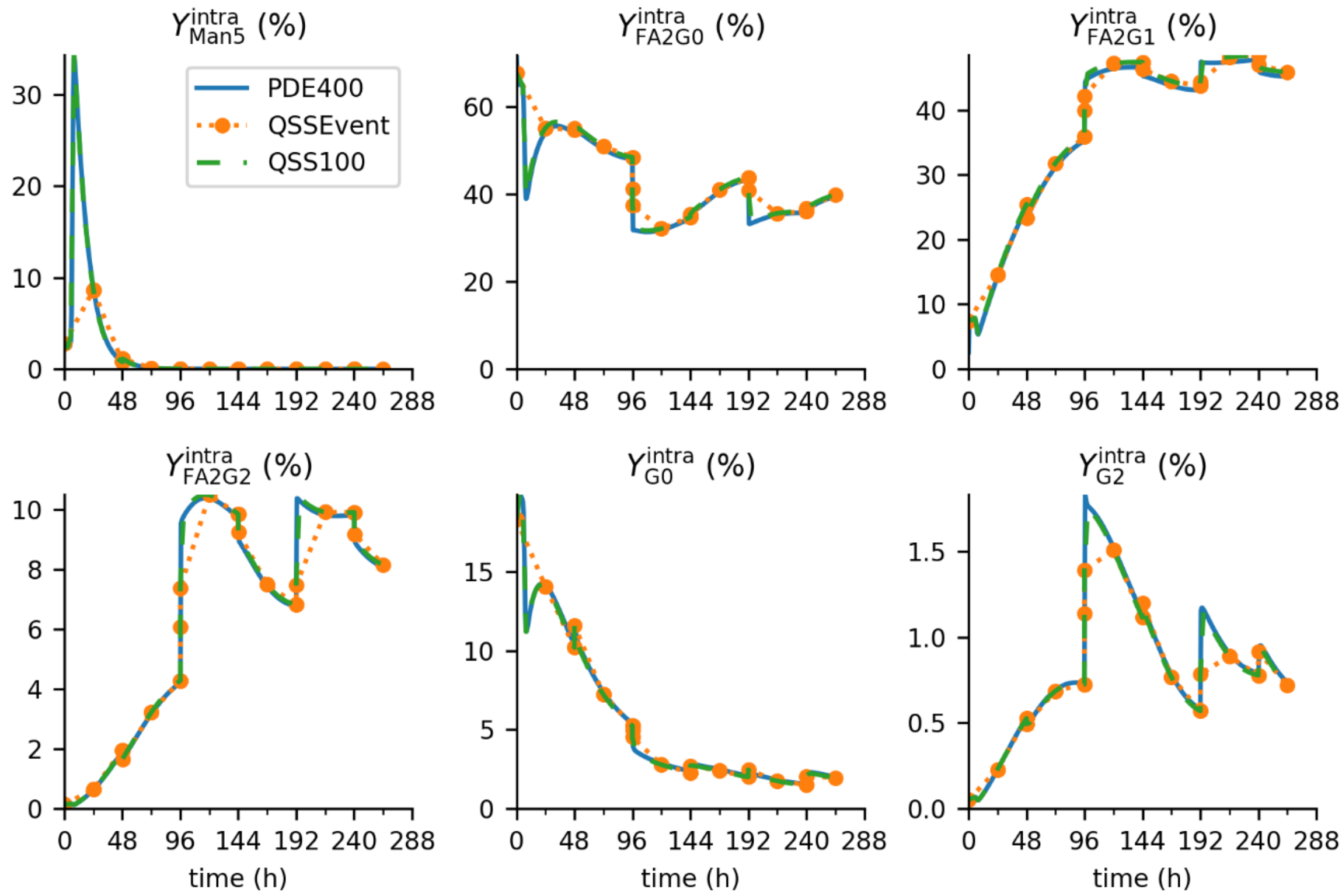

Figure 8. The trajectories of intracellular glycan compositions from PDE and QSS simulations.

Several conclusions can be drawn from Fig. 8. First, the QSS simulations at the time set $\mathcal{T}$ are accurate since all the circles representing the QSSEvent simulation results fall on the solid line denoting the PDE400 simulation. However, the dotted line has an evident mismatch from the solid line on the first day and the day after feeding galactose ($t = 96$ and $t = 192$ h) when constructing the trajectories of $Y_i^{\text{intra}}$ by connecting the circles, i.e., linear interpolation for the results at $\mathcal{T}_{\text{event}}$. This occurs because $\mathcal{T}_{\text{event}}$ is too sparse to capture the rapid changes in intracellular glycans during the time intervals where NSD concentrations continue to change quickly for 0.5 to several hours after the event time points. This is due to the mid-fast time scale of the NSD submodel, as mentioned in Section 3.1. However, with a denser time grid, the QSS100 simulation generates $Y_i^{\text{intra}}$ trajectories matching very well with those from the PDE400 simulation, as shown by the overlap between the solid line and the dashed line in Fig. 8. This capturing of the whole trajectories of intracellular glycan profiles is important as they will be used to compute extracellular glycan compositions.

Fig. 9 compares the errors of two QSS simulations with PDE simulations using different number of discretization points. Here, PDE50 and PDE100 refer to PDE simulations using 50 and 100 spatial discretization points, respectively. We assume the results from the PDE400 simulation as ground truth, and the absolute errors of the aforementioned methods are calculated by comparing them with these results.

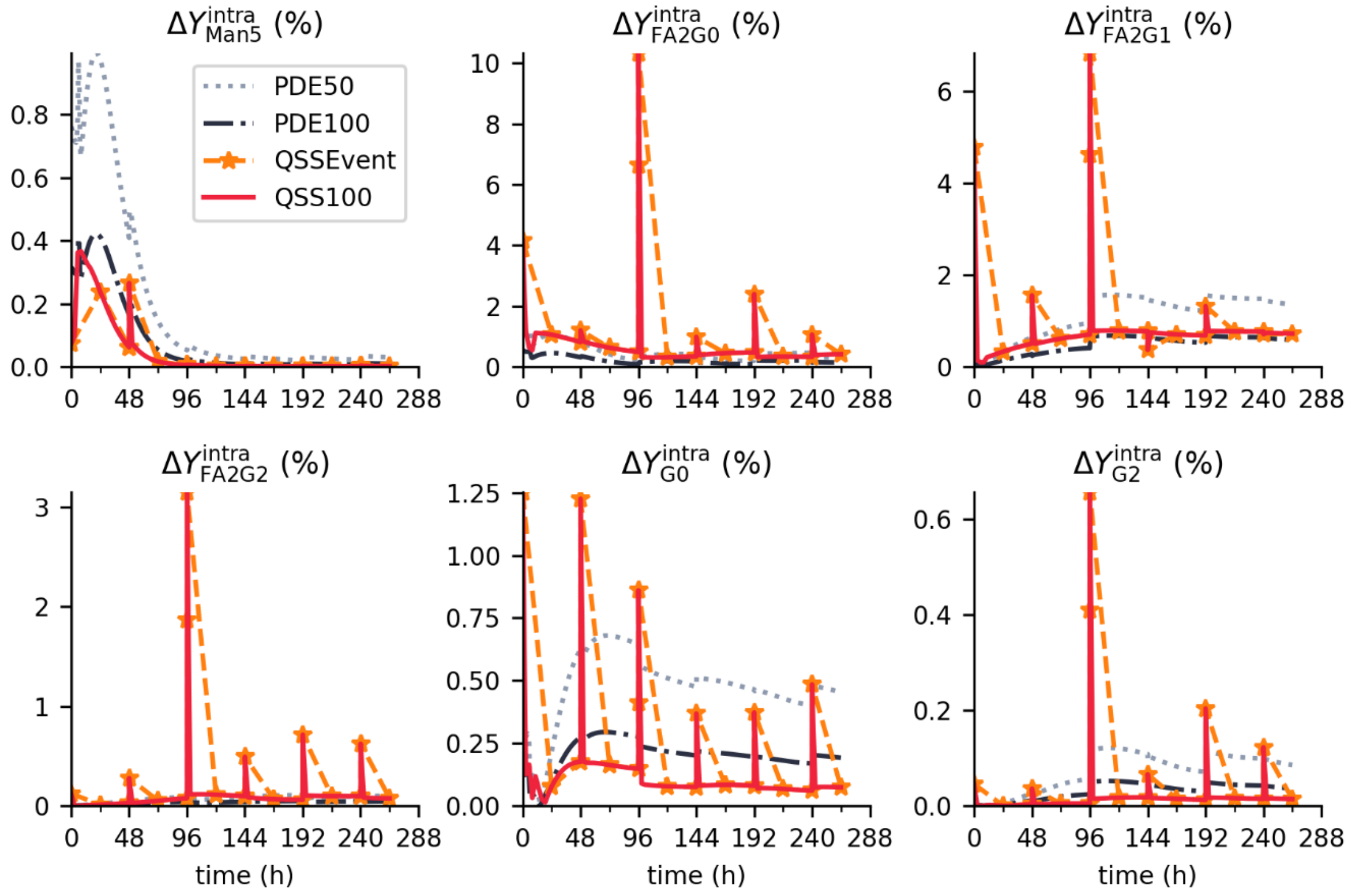

Figure 9. Computation errors of intracellular glycan compositions from PDE and QSS simulations.

As seen from Fig. 9, the QSS simulations usually have larger errors than all the PDE simulations during the initial $\tau$ and around the feeding time points, such as an error of 10% for $Y^{\text{intra}}_{\text{FA2G0}}$ around $t = 96$ h. This results from the sudden change in NSDs at those time points, causing an immediate response of $Y_i^{\text{intra}}$ in QSS simulations due to the disregard of the time delay ($20 \text{ min} < \tau < 30 \text{ min}$) in the steady-state Golgi model. However, in the PDE simulations, the time delay of the PFR model is considered. Nevertheless, such evident errors only appear in a short time interval (around $\tau$), so the accumulated errors in the computation of extracellular glycans is still small when the time grid is dense enough, e.g., in QSS100 simulation. On the other hand, in nearly all the other time points, the errors of $Y_i^{\text{intra}}$ in the QSS100 simulation are evidently smaller than those from the PDE50 simulation and close to those from the PDE100 simulation, as shown in Fig. 9. At those time points, the errors are less than 1% for FA2G0 and FA2G1, while they are less than 0.3% for the other glycans.

Overall, the above discussion shows two more error sources for the trajectories from QSS simulations: one is the interpolation between neighboring time points, and the other is the neglect of the time delay of the PFR model. The former can be addressed by adding more time points in $\mathcal{T}$, while the latter is negligible due to the short delay (less than 30 minutes).

**3.4 Simulation of the extracellular glycan trajectories**

The Golgi model predicts the intracellular glycan profiles at time points $\mathcal{T}$. To get the extracellular glycanforms, we need to integrate the cell culture model Eqs. (6) – (7), which requires the dynamic trajectories of $Y_i^{\text{intra}}$. Here, we treat $Y_i^{\text{intra}}$ as a time-variant parameter with a value $\bar{Y}_{i,k}^{\text{intra}}$ at $t_k \in \mathcal{T}_p = \mathcal{T} \backslash \{t_K\}$, which is computed from

$$\bar{Y}_{i,k}^{\text{intra}} = \frac{Y_{i,k}^{\text{intra}} + Y_{i,k+1}^{\text{intra}}}{2}, \qquad k \in \mathcal{T}_p \tag{13}$$

The computation of the extracellular glycans from intracellular glycans and cell culture model is shown in Fig. 10.

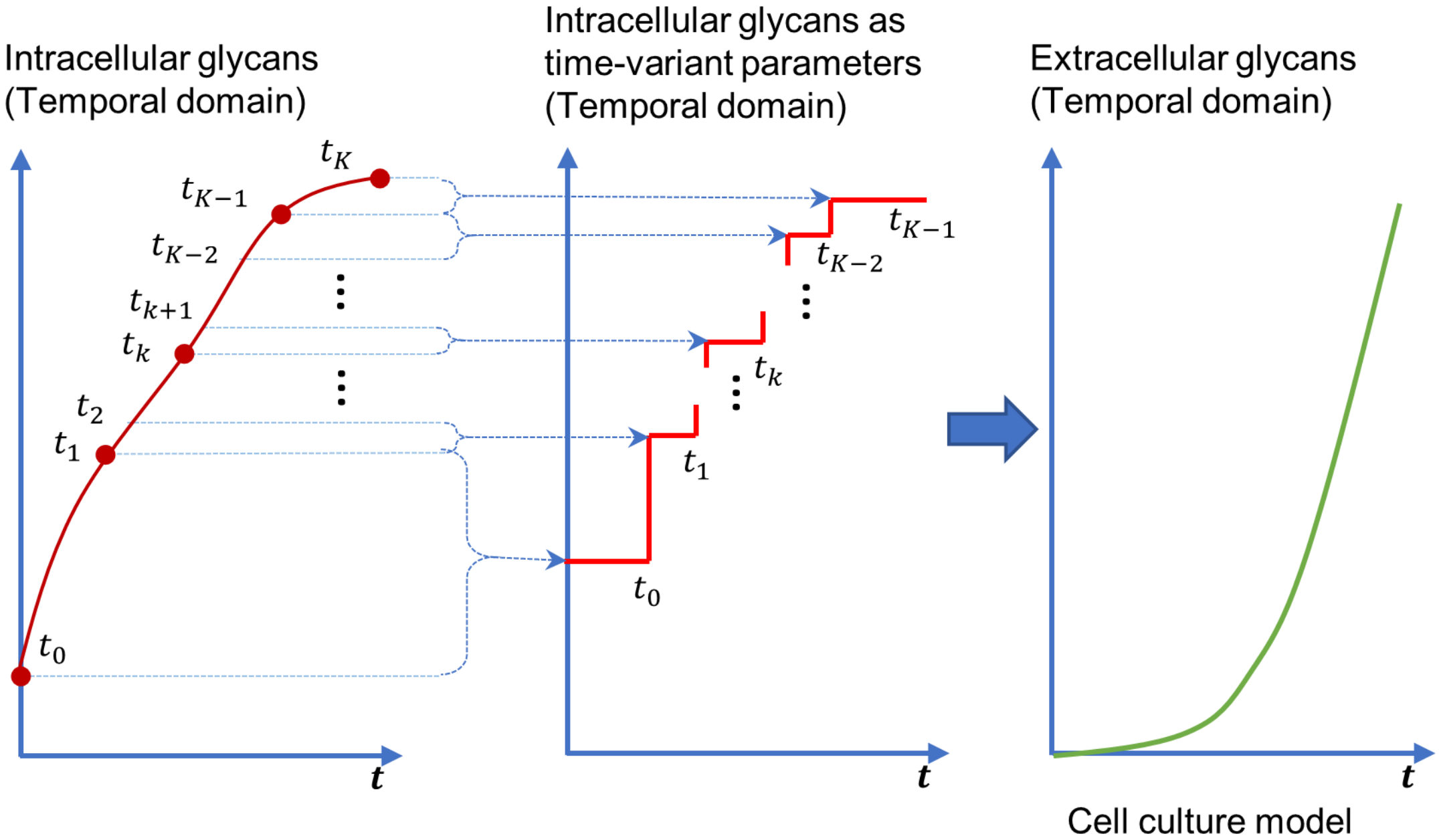

Figure 10. Computation of the extracellular glycan trajectories from intracellular glycan trajectories by treating intracellular glycans as time-variant parameters.

Although Eq. (13) is just a simple average of $Y_i^{\text{intra}}$ at two adjacent time points, it has been shown to be sufficiently accurate for approximating $Y_i^{\text{intra}}$ trajectories and for use in Eq. (6). This is demonstrated in Fig. 11, which compares extracellular glycan compositions from different simulation methods. Therefore, we will not bother with more sophisticated interpolation alternatives. As shown in Fig. 11, the extracellular glycan composition trajectories from the QSS100 simulation match those from PDE400 simulation very well. Moreover, the QSSEvent simulation in general has larger errors than the QSS100 simulation, especially in the first two days. Nevertheless, after two days, the QSSEvent simulation also generates good approximation to the PDE400 simulation results.

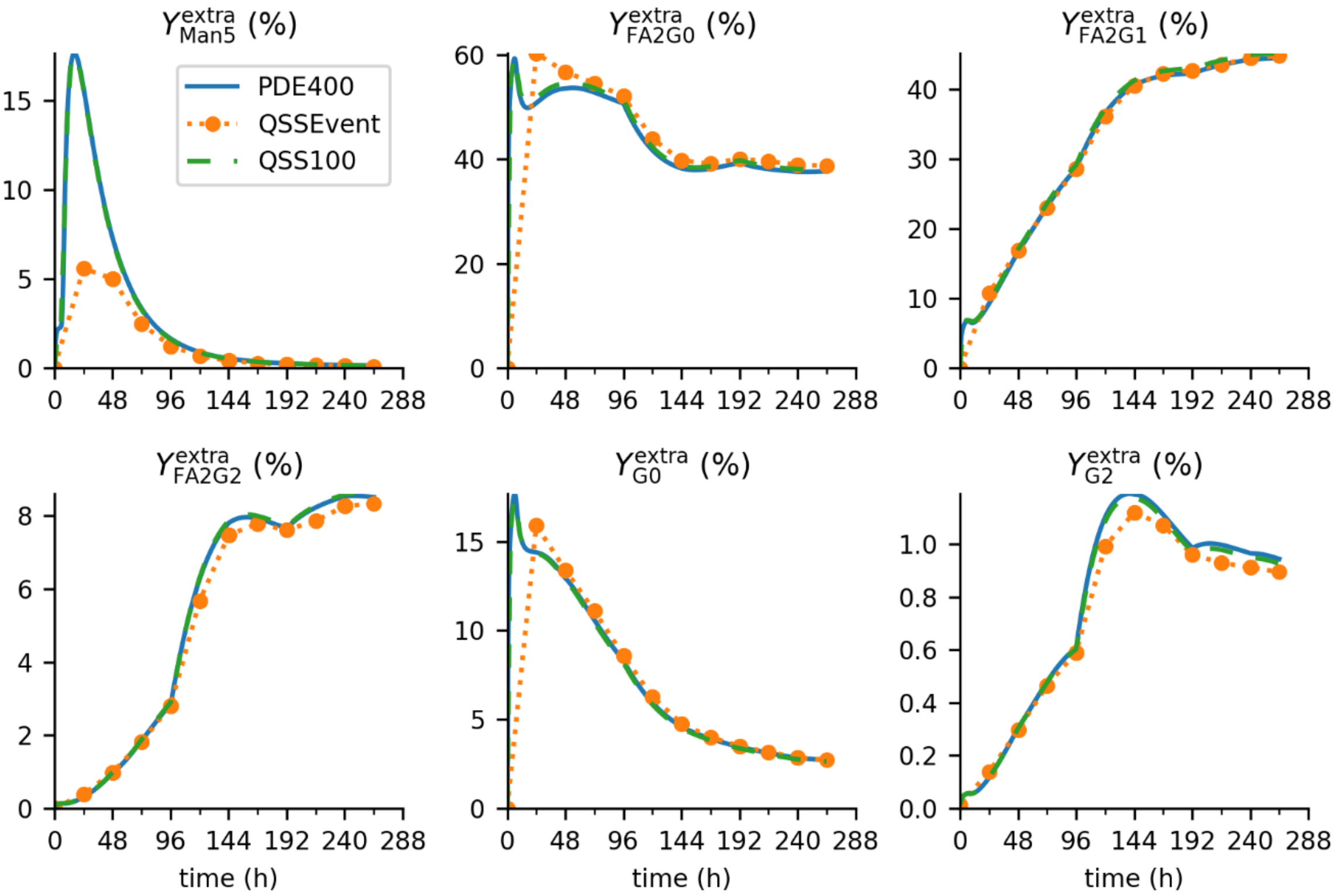

Figure 11. The trajectories of extracellular glycan compositions from PDE and QSS simulations.

The computational errors from different simulation methods are compared in Fig. 12. The errors of the QSS100 simulation are smaller than 1% for all the glycans throughout the cell culture period. Notably, at the end of the cultivation, the absolute errors are less than 0.7%. The largest error occurs in $Y^{\mathrm{extra}}_{\mathrm{FA2G1}}$, corresponding to a relative error of 1.6%. The errors of the QSS100 simulation are much smaller than those in the PDE50 simulation for all glycans except FA2G0. While the QSS100 simulation generates slightly larger error than the PDE100 simulation for major components (FA2G0 and FA2G1), it predicts the minor components (Man5, G0 and G2) with approximately 50% less errors.

Fig. 12 also shows that, in most cases, the QSSEvent simulation has larger errors than the other three simulations during the first two days. This is due to its sparse time grid for the Golgi simulations and the dramatic change in intracellular glycan profiles on the first day as shown in Fig. 8. However, after two days, its errors decrease to a level similar to the other simulations, with errors at most 2% for $Y^{\mathrm{extra}}_{\mathrm{FA2G0}}$ and less than 1% for the other glycans. For the fed-batch operation of mAb production using CHO cells, the glycoproteins produced in the first two days are usually not important. During this period, mAbs are far from being harvested and their concentrations are too low to be measured for parameter estimation and feedback control. Therefore, QSSEvent simulation can be considered for model-based analysis, optimization and control of the fed-batch bioreactor if QSS100 simulation is not fast enough.

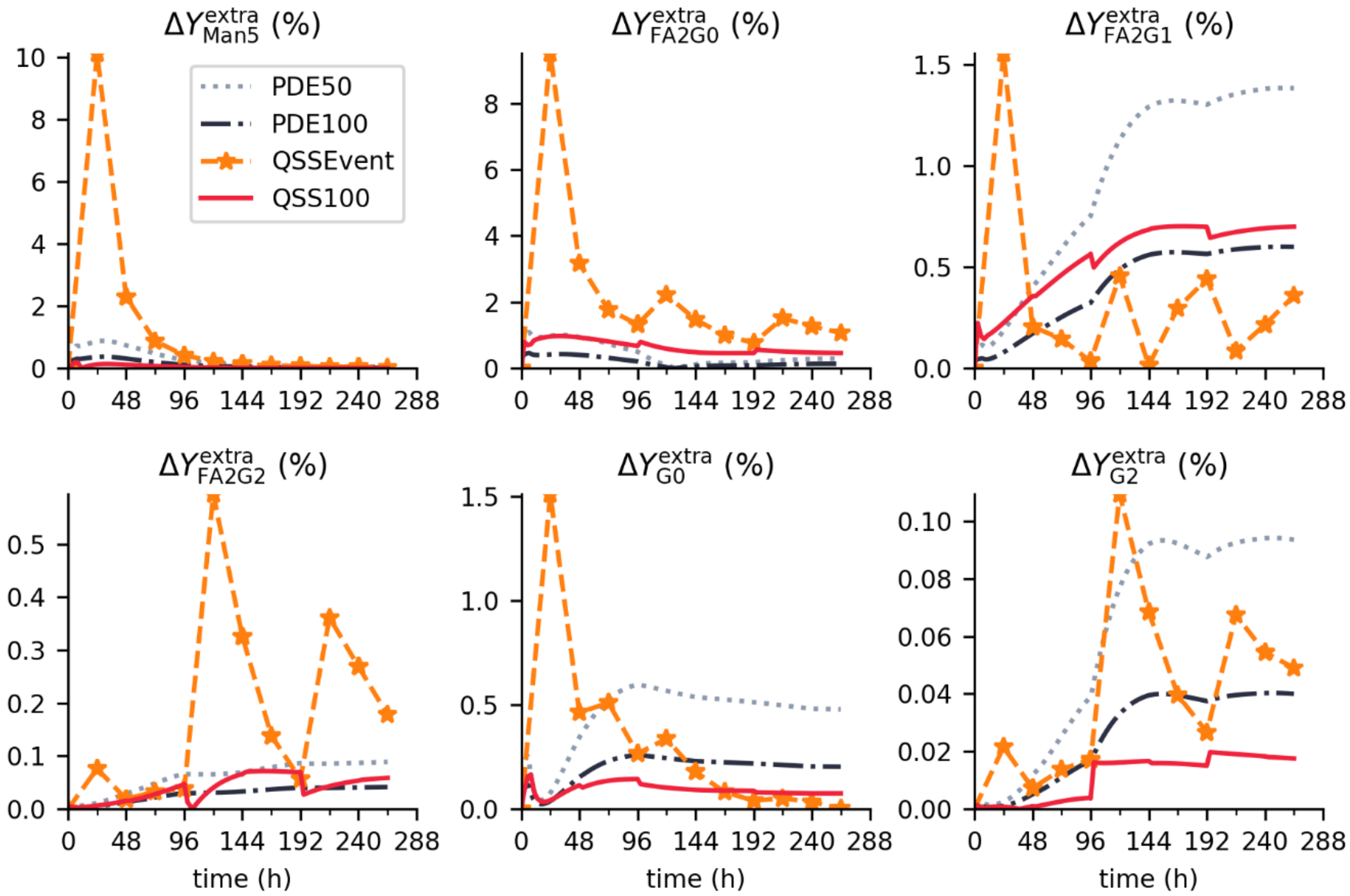

Figure 12. Computation errors of extracellular glycan compositions from PDE and QSS simulations.

### 3.5 Summary of the QSS simulation method

With above assumptions and analysis, the complete QSS simulation algorithm can be described as Algorithm 1.

**Algorithm 1**

Step 1: Set $\mathcal{T}$ and simulate the cell culture and NSD submodels with $\bar{Y}_i^{\text{intra}} = 0$; Obtain environment variables env at time points $\mathcal{T}$;

Step 2: Simulate the steady-state Golgi model at time points $\mathcal{T}$ to get the intracellular glycan profiles $Y_i^{\text{intra}}$ at the corresponding time points;

Step 3: Calculate $\bar{Y}_i^{\text{intra}}$ using Eq. (13), then integrate the cell culture submodel to obtain the extracellular glycan profiles $Y_i^{\text{extra}}$.

### 3.6 Influence of initial oligosaccharide concentrations on Golgi submodel simulations

Measuring the initial oligosaccharides concentrations inside the Golgi is difficult, so assumed initial values are typically used for simulating the Golgi submodel. However, since the residence time of oligosaccharides in the Golgi apparatus is short (less than half an hour), these initial oligosaccharides are quickly flushed out, having a negligible influence on the computation of subsequent intracellular glycan compositions. Additionally, the viable cell density is still very low at the beginning, so the secreted glycans during the short initial time interval has a minimal effect on the extracellular glycan compositions after one day. Therefore, using inaccurate initial values for the Golgi submodel is not a significant concern.

## 4. Speeding Up the QSS Simulation and Sensitivity Computation

Since our ultimate aim is to apply the glycosylation model for optimization and optimal control, both QSS and MOL methods for PDE simulations are implemented within the numerical optimal control framework, CasADi [27]. The DAE integrator IDAS is used for dynamic simulation [28]. The CasADi software and IDAS solver are also used to obtain the results in Section 3. Sensitivity equations are solved to obtain accurate derivatives, which are critical for efficient derivative-based optimization algorithms [29]. The forward sensitivity method is used here, as it has been found to have better convergence than the adjoint method for the multiscale glycosylation model. Two approaches for implementing the forward sensitivity method are compared. The first method involves calculating the sensitivity of dependent variables with respective to (w.r.t.) each input variable one by one and then concatenating all the sensitivities to obtain the entire Jacobian. The second method is to get the Jacobian of dependent variables w.r.t. all the input variables simultaneously by solving a much larger sensitivity equation system. These two methods are referred to as sensitivity_1by1 and sensitivity_simultaneous, respectively. Although the latter is faster, the former is less memory intensive, making it useful when the computer's memory is a bottleneck, as will be demonstrated sequentially. In this demonstration, the sensitivity computation involves the Jacobian of six extracellular glycan percentages w.r.t. 20 model parameters listed in Table S9. Note that the computation time for sensitivity computation usually increases with the number of independent variables for the forward sensitivity method. In this section, we fix the number to be 20, and it is enough to compare the performance of different algorithms.

All computations were conducted on a computer running the Windows 11 operating system, featuring a 12th Generation Intel® Core™ i7-12700H CPU with 14 physical cores and 20 logical cores, operating at 2.30 GHz, and equipped with 16 GB of RAM.

**4.1 Parallel computing**

In the QSS method, the steady-state Golgi submodel simulations at different time points are independent and can be conducted in parallel to save computational time. At each time point, the steady-state Golgi simulation involves solving a set of DAEs in which the "time" dimension corresponds to the position along the Golgi cisternae rather than chronological time. The only distinction among simulations is the vector of environmental variables env (Eq. 10), which are obtained from the simulation results of cell culture and NSD submodels. CasADi's map(•) function then provides a convenient mechanism for multithreaded execution of these independent steady-state Golgi simulations. Fig. 13 shows the computation time of QSS100 simulations using different numbers of threads.

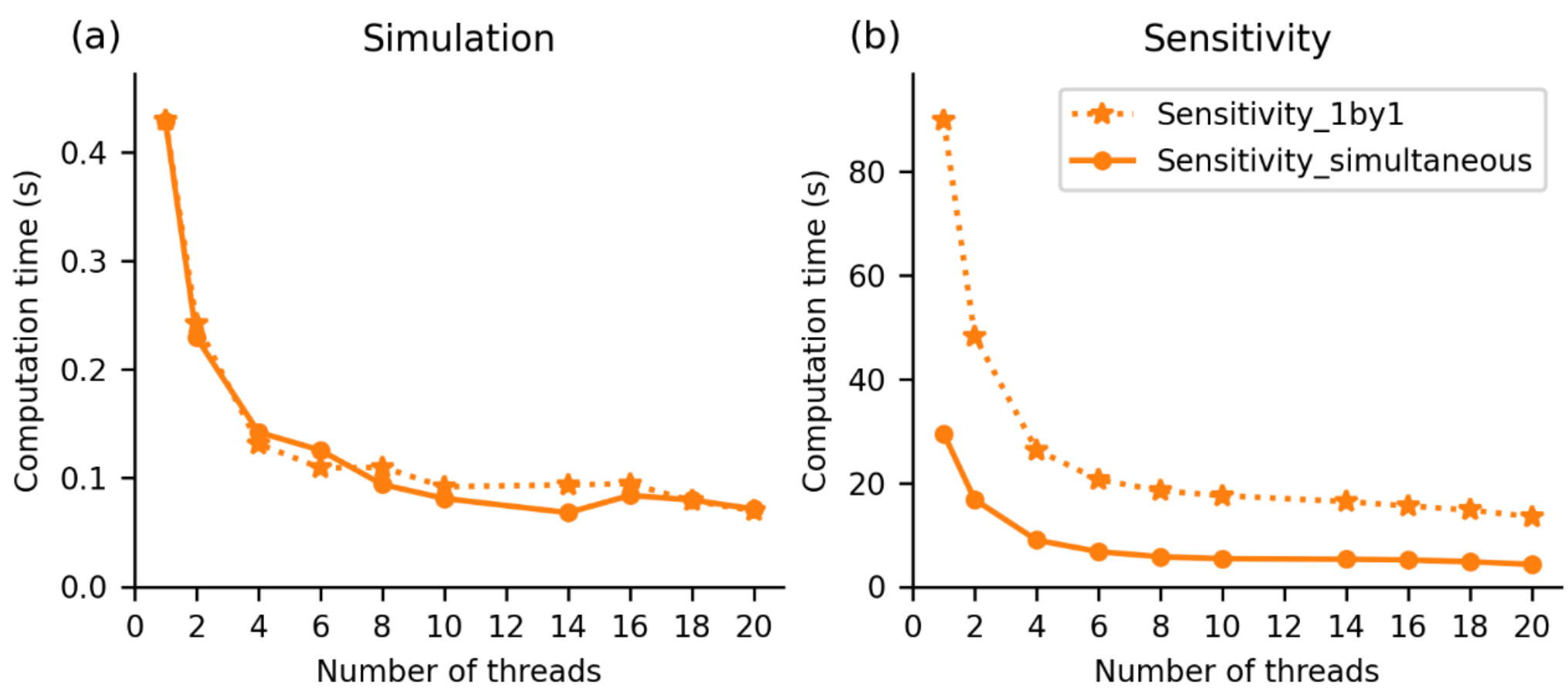

Figure 13. The computation time of (a) simulation and (b) sensitivity evaluation for QSS100 simulations using different numbers of threads.

As shown in Fig. 13, using more threads significantly improves the computational speed for both dynamic simulation and sensitivity computation. However, the benefit of additional threads diminishes as the number of threads $n_{\text{thread}}$ increases, particularly beyond $n_{\text{thread}} > 6$. The maximum time savings are around sixfold for simulation and sevenfold for sensitivity computation. Additionally, the simultaneous sensitivity evaluation method is approximately three times faster than the one-by-one sensitivity evaluation method as shown in Fig. 13b.

**4.2 Nonuniform time points allocation for the QSS simulation**

As seen in the last subsection, although we can use parallel computing to run multiple Golgi simulations simultaneously, there are usually not enough cores to run all the simulations at the same time on most personal computers. Therefore, it is still important to try to reduce the number of allocated time points for the Golgi simulations.

In the QSSEvent simulation, there are evident errors in the interpolated intracellular glycans trajectories during the first day and after feeding events. This is due to the rapid changes in NSD concentrations within half an hour to several hours after $t = 0$ and event time points, leading to dramatic changes in intracellular glycans during these periods. The QSS100 simulation captures these changes by adding 100 uniformly distributed time points, but a more parsimonious approach can be taken by selecting time points based on the time scales of the NSD submodel. As mentioned in Section 3.1, it takes around 6 hours for the concentration of UDPGlcNAc to reach to a quasi-steady state on the first day, and less than one hour for the NSD concentrations to stabilize after events. Therefore, in the QSS simulation, we can allocate more time points during these periods. The proposed strategy is to add nine more time points uniformly between 0 and 20 h and one point two hours after each event. Using two hours instead of one ensures that we do not miss the periods of rapid NSD concentration changes. The nonuniform time grid, denoted as $\mathcal{T}_{\text{nonuniform}}$, is shown in Fig. 14.

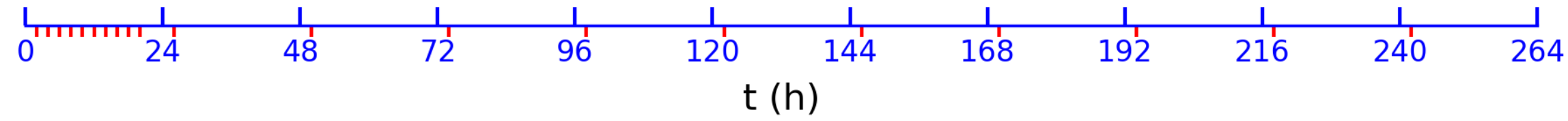

Figure 14. Nonuniform time grid for QSS simulation. The long ticks denote the time points in $\mathcal{T}_{\text{event}}$, while the short ticks denote the time points newly added in $\mathcal{T}_{\text{nonuniform}}$.

The simulation using the nonuniform time grid is denoted as QSSNonuniform and uses 48 time points for the experiment 10G. Fig. 15 presents that using half the number of time points in QSSNonuniform reduces the computation time by about a factor of 2 and up to 2.5 times compared to the QSS100 simulation.

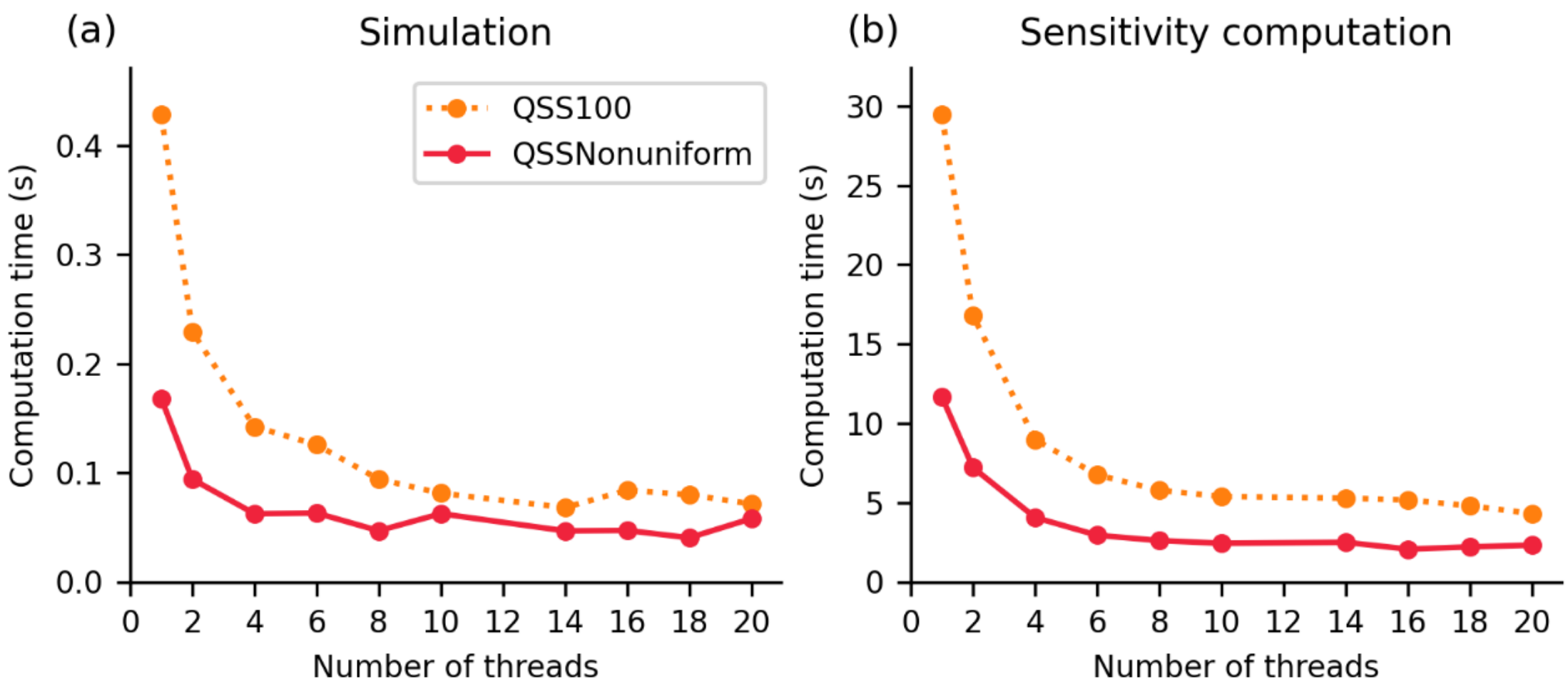

Figure 15. The computation time of (a) simulation and (b) sensitivity evaluation for QSS100 and QSSNonuniform simulations using different numbers of threads.

It is possible to reduce the number of discretization points in the PDE simulation without compromising accuracy by allocating discretization points more carefully. However, the optimal allocation in the spatial domain can vary significantly with different environment variables, making it less suitable for optimization problems. In contrast, the proposed time point allocation scheme in the QSS simulation remains constant under different simulation scenarios because the time scale for the NSD submodel after excitation is consistent. This is evidenced by the fact that the QSSNonuniform and QSS100 simulations produce nearly identical trajectories and similar errors for all five experiments, as shown in Figs. 16-17 and the Figs. S43-S47.

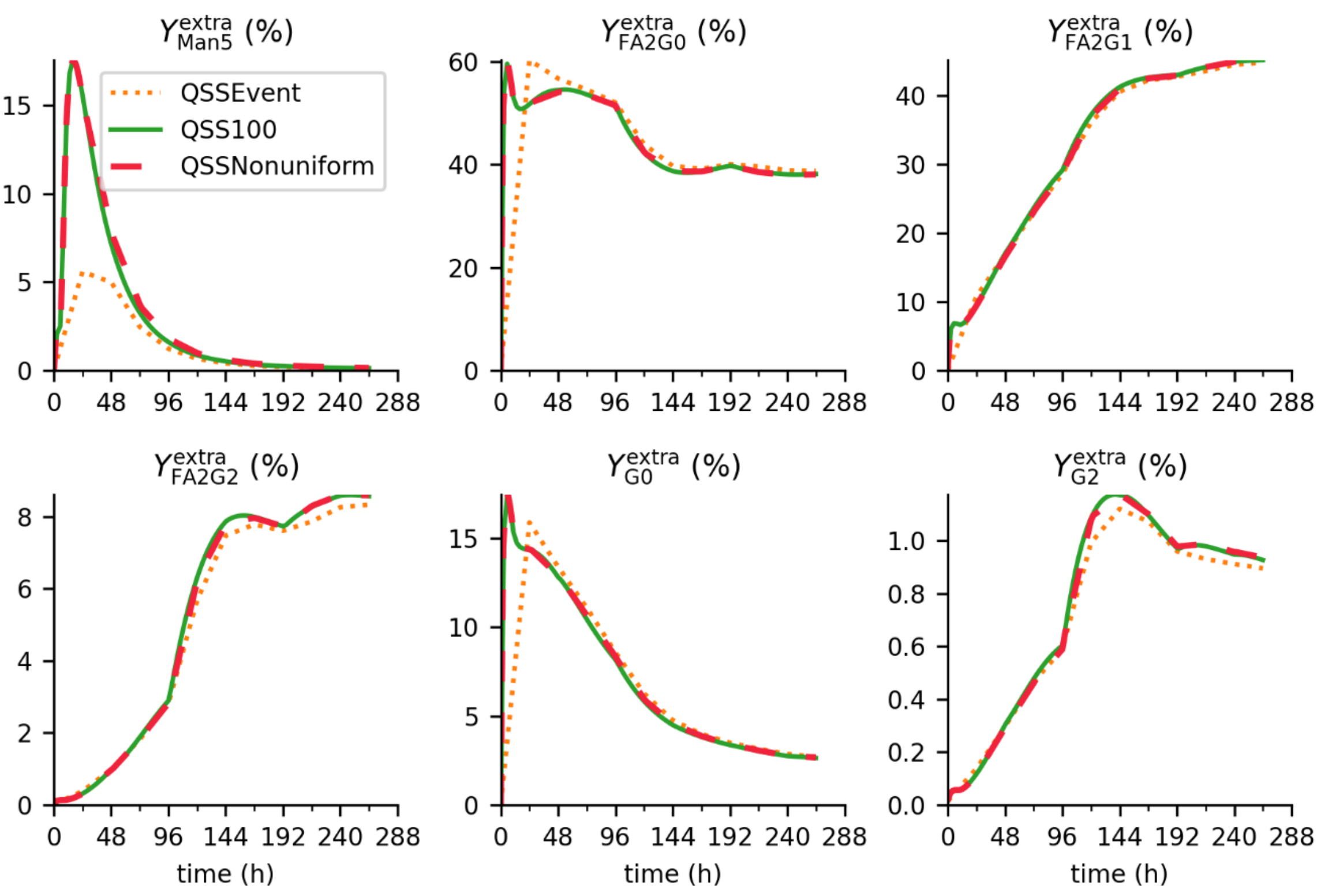

Figure 16. The trajectories of extracellular glycan compositions from QSS simulations with different time grids.

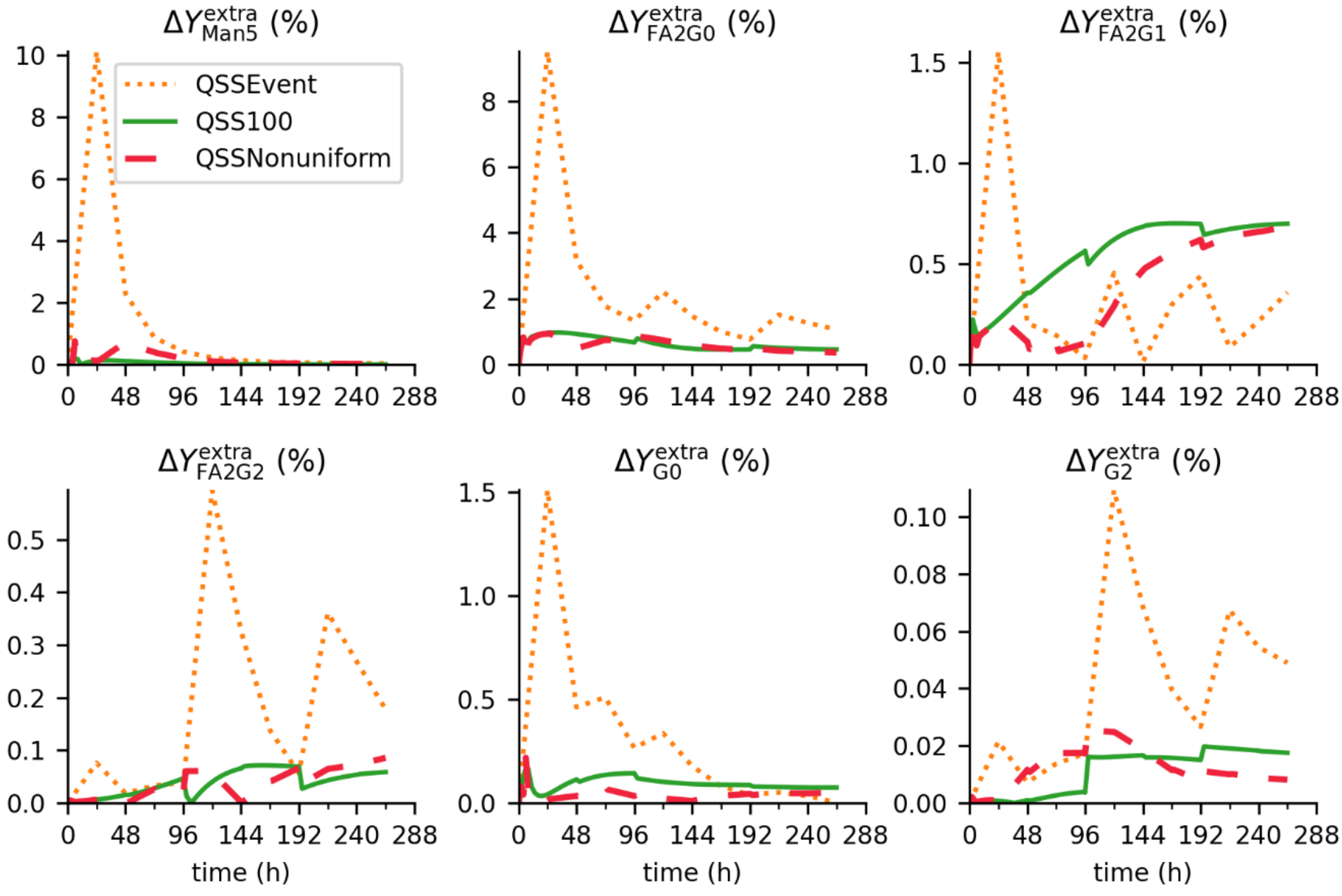

Figure 17. Computation errors of extracellular glycan compositions from QSS simulations with different time grids.

### 4.3 Computation efficiency comparison among PDE and QSS simulations

Since PDE100, QSS100 and QSSNonuniform simulations exhibit similar and satisfactory accuracy, as demonstrated in Section 3, we compare the computation time for these three simulations, as shown in Fig. 18.

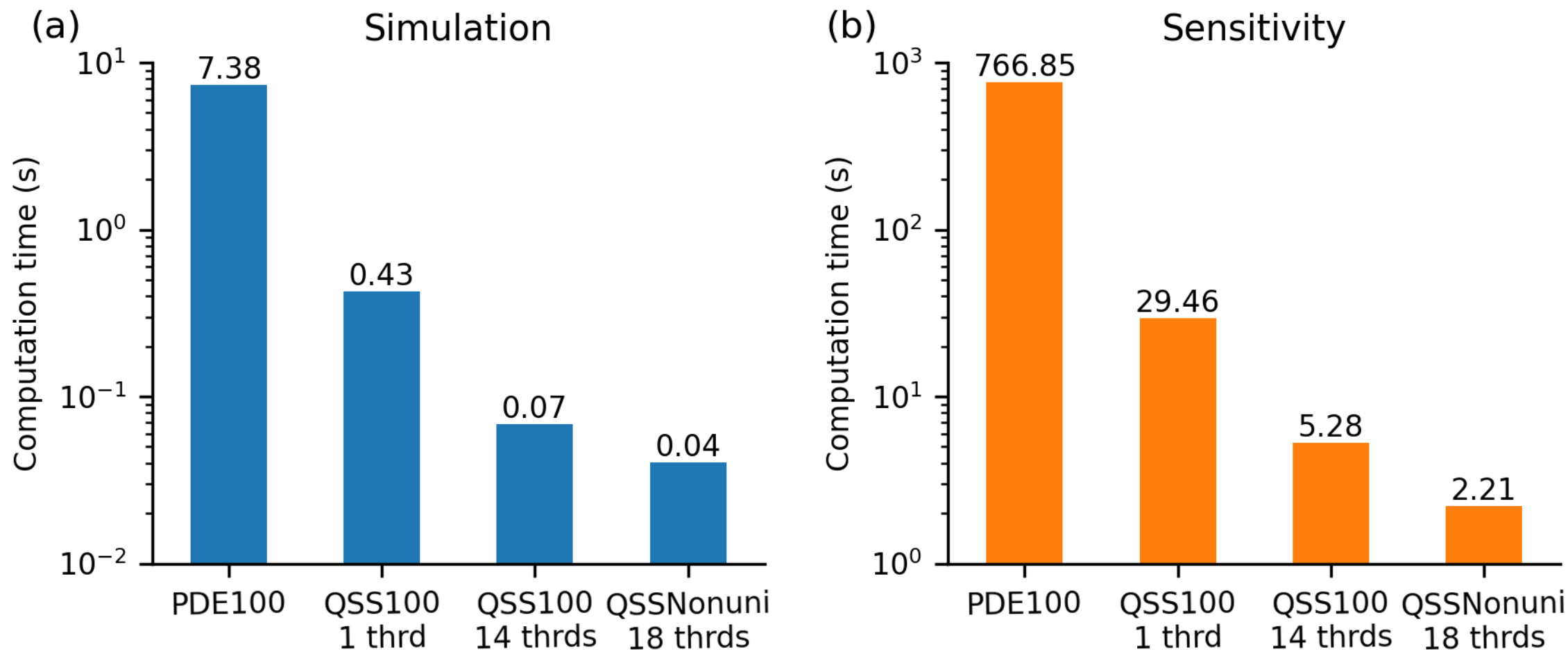

Figure 18. Computation time of (a) simulation and (b) sensitivity evaluation when using different simulation methods. QSS100 1 thrd: QSS100 simulation using 1 thread; QSS100 14 thrds: QSS100 simulation using 14 threads; QSSNonuni 18 thrds: QSSNonuniform simulation using 18 threads. PDE100 simulation cannot be parallelized.

As shown in Fig. 18a, the QSS100 simulation reduces computation time by around 17 times compared to the PDE100 simulation. Furthermore, the application of parallel computing to QSS100 simulation achieves a time saving of around 100 times. The combined use of parallel computing and the nonuniform time grid results in the greatest time saving compared to the PDE100 simulation, achieving approximately 180 times better performance.

Fig. 18b shows that it takes over 12 minutes to compute the sensitivity of six glycans w.r.t. 20 model parameters using the PDE100 simulation, even though the simulation itself only takes 7 seconds. This underscores the strong motivation to improve computation speed. The QSS100 simulation alone, without parallel computing, significantly reduces the computation time to less than half a minute, achieving a time saving of about 26 times compared to the PDE100 simulation. Using 14 threads, the QSS100 simulation and the QSSNonuniform simulation further reduce the computational time to 5 seconds and 2 seconds, respectively, demonstrating time savings of 145 times and 348 times.

From the above discussion, it is evident that QSS simulations achieve greater time savings for sensitivity evaluation than for simulation. This is because the sensitivity_1by1 method must be used for sensitivity evaluation based on the PDE100 simulation to avoid the memory limitations encountered by the faster but more memory-intensive sensitivity_simultaneous method. In contrast, the efficient simultaneous sensitivity method can be applied to the QSS simulations without encountering the memory issues.

## 5. Optimization Problems Using QSS Simulations

In Section 4, we validated the proposed QSS simulation algorithm solely for simulations and sensitivity computation. In this section, we will demonstrate the accuracy and efficiency of the QSS algorithm in two important optimization problems: parameter estimation and dynamic optimization.

The control vector parameterization (CVP) method is used for both optimization problems due to its good convergence [30]. A robust and efficient sequential quadratic programming (SQP) implementation called PySQP, based on Ma et al. [31] with watchdog technique [32] improvement, is used to drive the CVP method. The derivatives used in the SQP algorithm are obtained by solving sensitivity equations as mentioned in Section 4. Since dynamic simulations and sensitivity computations are conducted repeatedly in the CVP method, these components dominate the total computational time. This motivates us to develop the efficient and accurate QSS simulation method. In this section, we use the QSS100 simulation with multiple threads for the optimization, as it has been proven to be sufficiently fast.

### 5.1 Parameter estimation

Parameter estimation is typically the first step in model-based analysis, obtaining the parameters used for simulations and optimization. Moreover, parameter estimation is a crucial step in adaptive control. Therefore, ensuring both the accuracy and efficiency of parameter estimation is essential.

In this case study, we conduct parameter estimation using data from five experiments by Kotidis et al. [9] due to the richness of the dataset. For the parameters in the cell culture and NSD submodels, as well as the enzyme distribution parameters in the Golgi submodel, we directly adopt those from Kotidis et al. [9]. However, we estimate all the dissociation parameters and kinetic constants in the Golgi submodel, resulting in 20 parameters to be estimated, as shown in Table S9. The maximum a posterior (MAP) method is used for parameter estimation because it can utilize existing literature data and mitigate the unidentifiability issue [33]. The MAP estimation can be formulated as the minimization,

$$\min_{\beta} \left[ \left(Y - \hat{Y}\right)^{\mathrm{T}} V_{\epsilon}^{-1} \left(Y - \hat{Y}\right) + (\beta - \mu)^{\mathrm{T}} V_{\mu}^{-1} (\beta - \mu) \right], \tag{14}$$

where $\beta$ and $\mu$ are the parameter estimates and prior parameters, respectively, $Y$ and $\hat{Y}$ are the experimental data and model predictions, respectively, and $V_{\epsilon}$ and $V_{\mu}$ are measurement covariance and prior parameter covariance, respectively. The prior parameter values and their variances are obtained from Kotidis et al. [9] and Karst et al. [22], as shown in Table S9. The initial values (also prior values), lower bounds, and upper bounds of the parameters are also provided in Table S9. The optimal objective function values and the computation times from estimations using QSS and PDE simulations are shown in Fig. 19. In Fig. 19a, the objective function values are obtained by substituting the parameters derived from different simulation method-based estimations into the PDE400 simulation for a fair comparison.

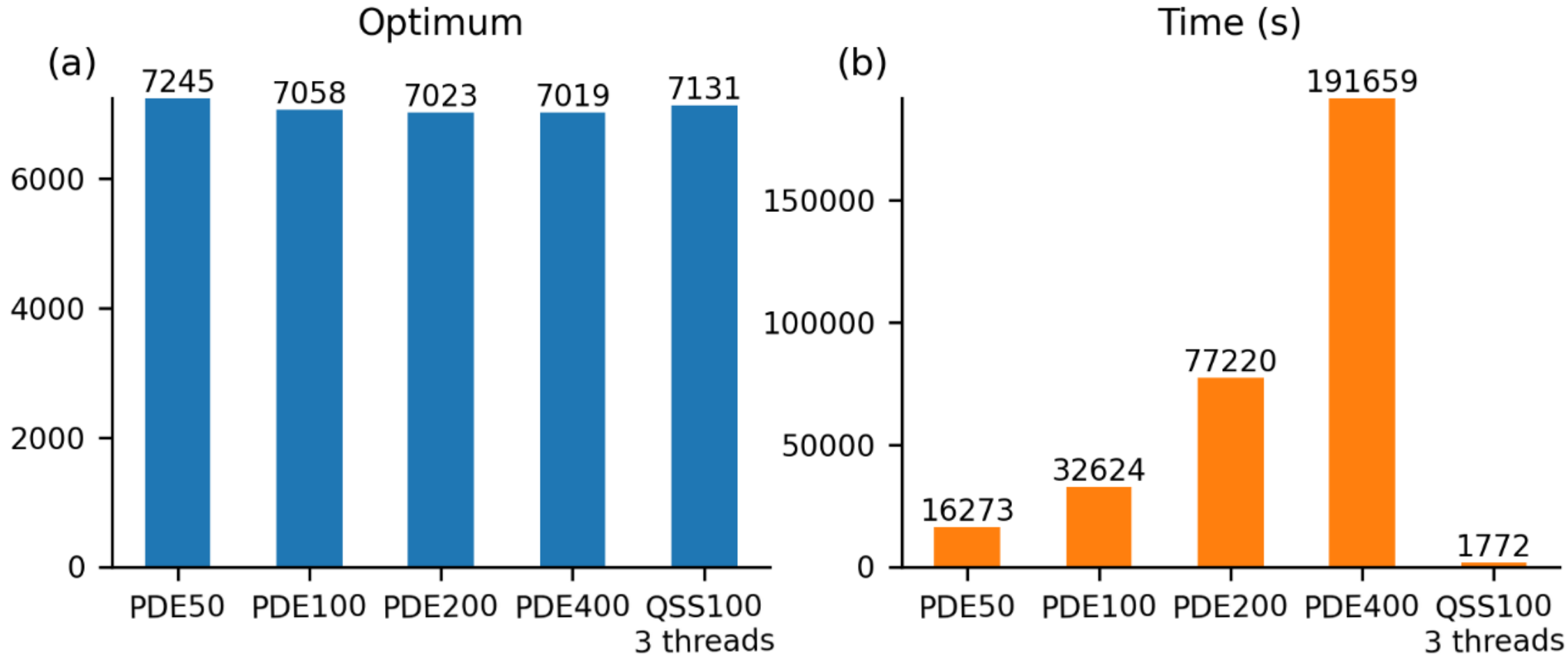

Figure 19. a) Optima and b) computational times of the parameter estimations using different simulations methods.

According to Fig. 19a, all methods generate similar optima. The largest minimum is obtained from the parameter estimation using the PDE50 simulation, which is relatively 3% larger than that from the PDE400 simulation. The QSS100 simulation-based estimation produces an optimum of 1.6% larger than the PDE400 simulation-based estimation. In terms of the total computational time, as shown in Fig. 19b, the QSS100 simulation-based estimation achieves the optimum in less than 0.5 hour, while the PDE simulation-based estimations require 4.5 to 53 hours, presenting a time saving of one to two orders of magnitudes. Note that the parallel simulations of the five experiments (requiring five cores) are applied to both PDE and QSS simulation-based parameter estimations. Therefore, the benefit of using parallel computing for each QSS simulation will not be evident on our personal computer for the current parameter estimation problem, considering that the improvement in computation speed is nonobvious after utilizing more than six threads, as mentioned in Section 4.1.

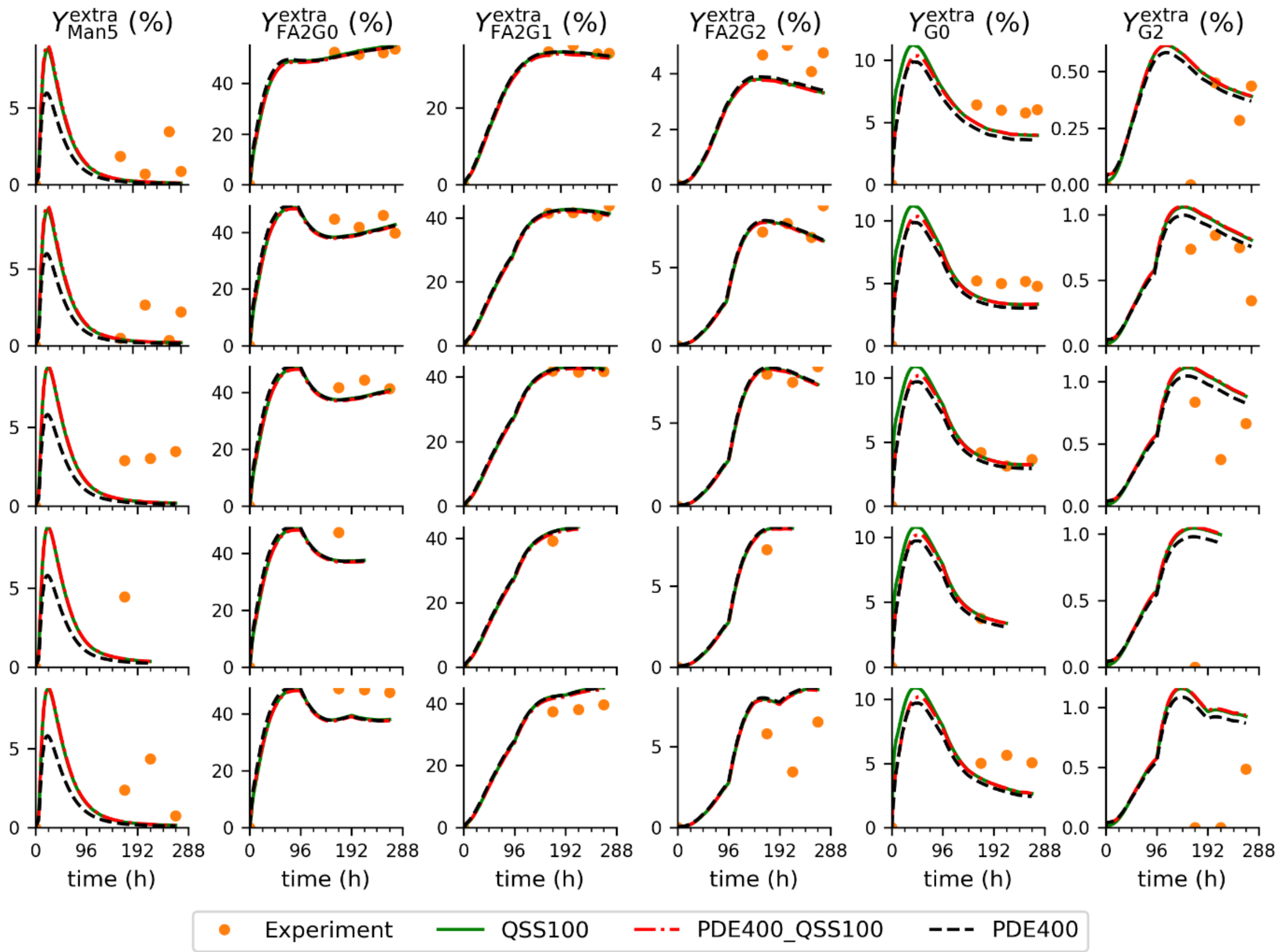

Figure 20. The comparison between experimental data and different simulations using estimated parameters. QSS100 denotes the QSS100 simulation results using parameters from QSS100 simulation-based estimation; PDE400_QSS100 denotes the PDE400 simulation using parameters from QSS100 simulation-based estimation; PDE400 denotes the PDE400 simulation using parameters from PDE400 simulation-based estimation.

Fig. 20 compares experimental data with the results from QSS100 simulation using its own estimated parameters (denoted as QSS100), PDE400 simulation using its own estimated parameters (denoted as PDE400), and PDE400 simulation using parameters from QSS100 simulation-based estimation (denoted as PDE400_QSS100). All methods generate good fits to the experimental data, especially for the major components. The results from QSS100 and PDE400_QSS100 simulations are nearly identical, except for the prediction of G0 in the first three days. Even for these values, the deviations between the two simulations are still less than 1% excluding the first several hours. The predicted trajectories from PDE400_QSS100 and PDE400 overlap very well for five out of the six glycans throughout the entire cell cultivation, with evident deviations in the Man5 percentage in the first 7 days. For that period, it is difficult to determine which fitting is more accurate due to the lack of experimental data. Overall, both the accuracy of the QSS simulation and the QSS simulation-based parameter estimation are satisfactory.

**5.2 Dynamic optimization with path constraints**

The calibrated model can be employed to enhance the performance of our bioreactor through model-based dynamic optimization. In this study, we focus on optimizing a fed-batch experiment conducted in a shaking flask. The flask runs for 12 days, with a pulse feed at the beginning of each day and a 10 mL sample taken at the end of each day (the total sampling time is 36 s). The product, mAb, is harvested at the end of the cultivation. The initial working volume of the shaking flask is 100 mL, with allowed variations between 75 and 150 mL during cultivation. The decision variables for the optimization problem are the feed flow rate ($F_{in}$), the concentration of galactose ($[\mathrm{Gal_{feed}}]$) in the feed stream, and the concentration of uridine ($[\mathrm{Urd_{feed}}]$) in the feed stream. To mimic the pulse feed, feeding is assumed to occur over 36 s. The concentrations of the other substrates are shown in Table S3.

The objective of the optimization is to maximize the concentration of the galactosylated mAb species, calculated using [9]

$$\mathrm{galactoylated_mAb} = [\mathrm{FA2G1}] + 2 \cdot [\mathrm{FA2G2}]. \quad (15)$$

Constraints are added to the optimization to ensure a feasible process. These constraints include a lower bound of cell viability (viability) and lower and upper bounds on the working volume ($V$, mL). Additionally, there are bound constraints for the operational variables $u := [F_{\mathrm{in}}, [\mathrm{Gal_{feed}}], [\mathrm{Urd_{feed}}]]$, i.e., $u \in [u_{\mathrm{lb}}, u_{\mathrm{ub}}]$. The constraints are detailed in the optimization,

$$\max_{u(t)} \mathrm{galactosylation} \quad (16)$$

$$\text{s.t.} \quad x(0) - x_0 = 0, \quad (17)$$

$$\mathrm{Glyco}(x(t), z(t), u(t), p) = 0, \quad (18)$$

$$\mathrm{viability}(t) \geq 60\%, \quad (19)$$

$$75 \leq V(t) \leq 150, \tag{20}$$

$$u_{\mathrm{lb}} \leq u(t) \leq u_{\mathrm{ub}}, \tag{21}$$

$$t \in [0, T],$$

where $\mathrm{Glyco}(\cdot)$ is the multiscale glycosylation model and all path constraints are enforced throughout the cell cultivation period. The initial values, lower bounds, and upper bounds of the decision variables are given in Table 2.

Table 2. Initial values, lower bounds, and upper bounds of the decision variables.

| Variable names | $F_{\mathrm{in}}$ ($\mathrm{L \cdot h^{-1}}$) | $[\mathrm{Gal_{feed}}]$ (mM) | $[\mathrm{Urd_{feed}}]$ (mM) |
|---|---|---|---|
| Initial values | 1 | 1 | 1 |
| Lower bounds | 0 | 0 | 0 |
| Upper bounds | 100 | 1000 | 1000 |

The dynamic optimization is solved by the CVP method. In total, there are 36 decision variables. When the sensitivity_simultaneous method was used to compute the derivatives in the PDE model-based dynamic optimization, CasADi terminated prematurely due to insufficient memory on the personal computer, which cannot accommodate the symbolic Jacobian. Therein, sensitivity_1by1 method or finite difference method has to be applied for the optimization. There is no memory issue for QSS simulation-based dynamic optimization, allowing the efficient simultaneous sensitivity computation method to be used.

Fig. 21 compares the optimal objective function values and computation times from different simulation methods, including the QSS100 simulation using 15 threads, PDE simulations using the finite difference method for derivatives, and PDE simulations using sensitivity_1by1 method for derivatives. As with the last case study, for a fair comparison, the objective function values in the figure are obtained by substituting the optimal decision variables from different optimizations into the PDE400 simulation. As seen in Fig. 21a, the PDE-based optimizations using the finite difference method either generate a solution with significantly lower galactosylation (around 12%) than the other optimizations or terminate at an infeasible solution. This highlights the importance of using sensitivity equations to provide accurate derivatives for optimization. PDE50 and PDE100 simulation-based optimizations using sensitivity_1by1 and QSS100 simulation-based optimization obtain very similar solutions, although the latter is 0.1% lower. In terms of computation time, QSS100 simulation-based optimization needs less than 15 minutes, while PDE50 and PDE100 simulation-based optimizations require around 20 hours and 28 hours, respectively, as shown in Fig. 21b. This demonstrates a speed improvement of around 80 times and 115 times speed improvement, respectively. Given the slow performance of PDE100 simulation-based optimization using sensitivity_1by1 and the satisfactory results of the QSS100 simulation-based optimization, we did not conduct optimization using more discretization points for the PDE simulation-based optimization.

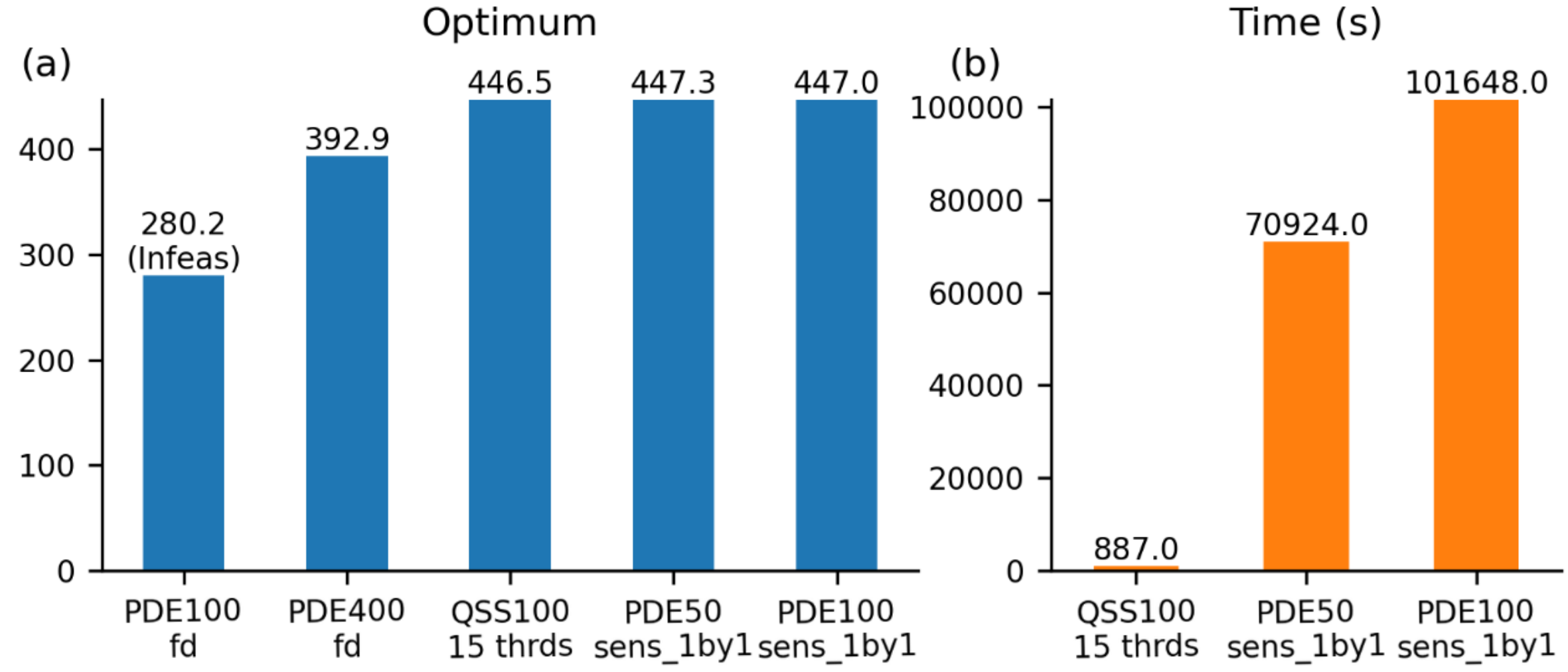

Figure 21. a) Optima and b) computation times of optimizations using QSS and PDE simulations. PDE100-fd, PDE400-fd: PDE simulation-based optimizations using the finite difference method for derivatives; QSS100-15 thrds: QSS100 simulation using 15 threads; PDE50-sens_1by1, PDE100-sens_1by1: PDE simulation-based optimizations using the sensitivity_1by1 method for derivatives.

After QSS100 simulation-based optimization, the concentration of galactosylated mAb increases to 446.5 $\mathrm{mg \cdot L^{-1}}$, which is 64% higher than the concentration (273.1 $\mathrm{mg \cdot L^{-1}}$) obtained in the experiment 10G, the highest among the five experiments used for model calibration. The optimal decision variables, and the corresponding optimal trajectories of some key state variables obtained from the QSS100 simulation and the PDE400 simulation (denoted as PDE400_QSS100) are shown in Fig. 22. According to Fig. 22a, nutrient supplements are added at four time points: the beginning of the first, second, third and eighth days. Note that there is no feeding when the feed flow rates are 0, even if the galactose and uridine concentrations are nonzero. All the feedings at the four time points contain galactose, but no uridine. This is because, according to Eqs. (S5) and (S7) and Grainger and James [15], uridine has a detrimental influence on cell viability, which can lead to lower cell density and unsatisfied viability constraints (Eq. 19). Conversely, galactose can increase UDP-Gal, which is required for galactosylation, and simultaneously benefit the viable cell density [15]. To demonstrate the importance of galactose on galactosylation, Fig. 22b also shows the state variable trajectories from the QSS100 simulation without galactose in the feed at $t = 0$, which presents a 47% decrease in the concentration of the galactosylated mAb. Finally, in Fig. 22b, the good match between the results from the QSS100 simulation and the PDE400_QSS100 simulation confirms the accuracy of the QSS100 simulation.

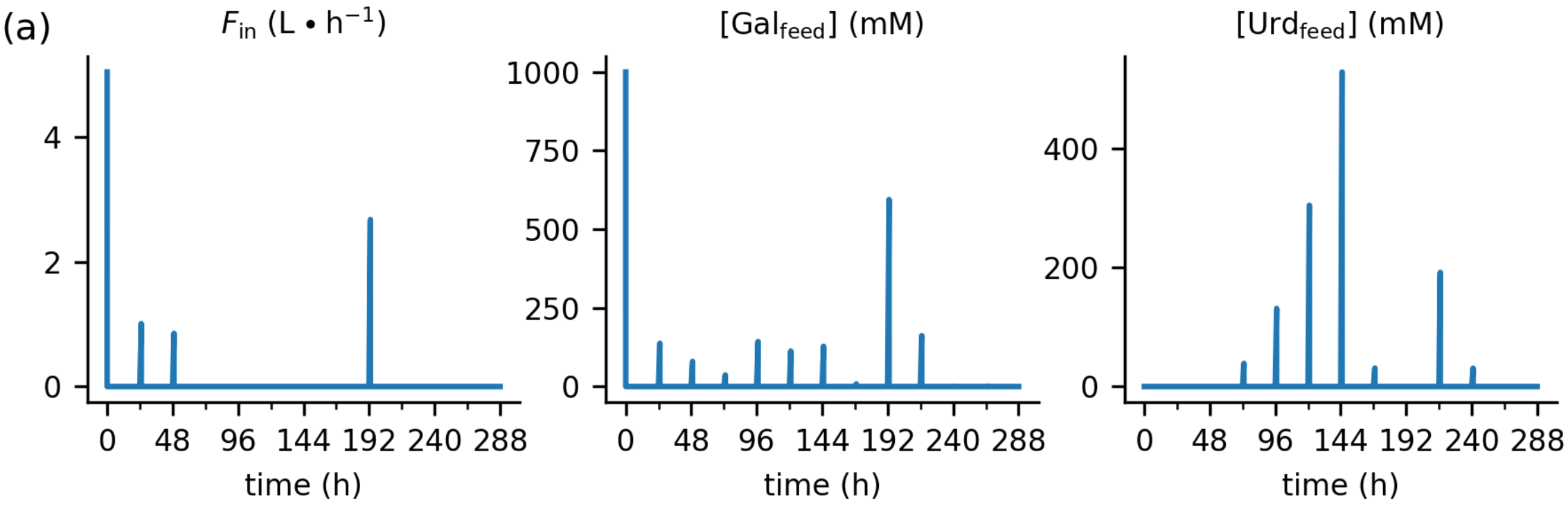

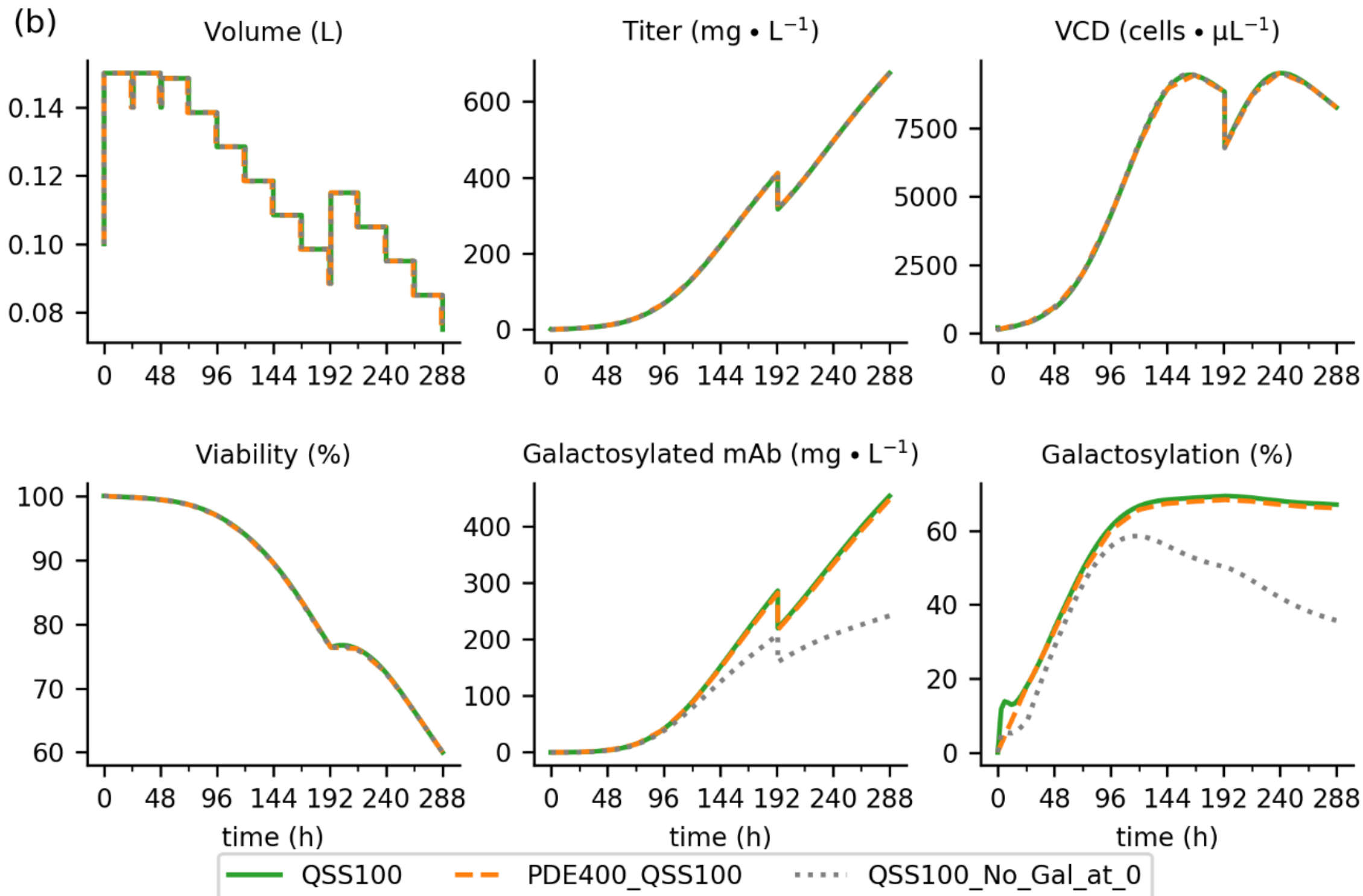

Figure 22. a) The optimal decision variable trajectories generated from the QSS100 simulation. b) The trajectories of key state variables when implementing the optimal decisions and comparison decision variables. QSS100: the optimal trajectories from the QSS100 simulation-based optimization; PDE400_QSS100: the trajectories from the PDE400 simulation using the optimal decision variables obtained from QSS100 simulation-based optimization; QSS100_No_Gal0: the trajectories from the QSS100 simulation using the optimal decision variables from QSS100 simulation-based optimization but no galactose feed at $t = 0$.

## 6 Conclusions

In this article, we propose a QSS simulation method for the efficient and accurate solution of the multiscale glycosylation model, which involves cell culture, intracellular NSD synthesis and Golgi glycosylation submodels. We introduce the QSS assumption for the Golgi submodel, which convert the PDAE model to a DAE model. Additionally, we assume negligible NSD consumption in glycosylation within the Golgi apparatus, decoupling the steady-state Golgi submodel in the spatial domain from the other submodels in the temporal domain, allowing sequential simulations. The intracellular glycan trajectories obtained from the Golgi submodel are then used in the cell culture submodel to compute the extracellular glycan compositions. Comparisons between QSS and PDE simulations for five experiments validate the assumptions and accuracy of the QSS simulation method.

To further improve the computational speed of the QSS simulation method, we apply parallel computing to the independent steady-state Golgi model simulations at different time points. Additionally, the Golgi simulations are conducted on a nonuniform time grid, reducing the number of time points while maintaining solution accuracy. The enhanced QSS method demonstrate a 180-fold improvement in simulation efficiency and a 348-fold improvement in sensitivity computation efficiency compared to the PDE simulation.

We also apply the QSS simulation method to two important optimization problems in the glycosylation process: parameter estimation and dynamic optimization. Case studies show that the QSS simulation-based optimizations are one to two orders of magnitude faster than PDE simulation-based optimization, with negligible solution degradation.

The proposed QSS simulation method provides a solid foundation for fast analysis and optimization of the mAb glycosylation process, enabling future work in moving horizon estimation, nonlinear model predictive control, and model-based experimental design.

**Acknowledgement**

This study was financially supported by a Project Award Agreement from the National Institute for Innovation in Manufacturing Biopharmaceuticals (NIIMBL), U.S., with financial assistance from awards 70NANB17H002 and 70NANB20H037 from the U.S. Department of Commerce, National Institute of Standards and Technology. We would also like to thank Professor Cleo Kontoravdi for answering our questions on the Golgi glycosylation model and for generously providing the experimental data in her paper.

**Nomenclature**

| Abbreviations | Definition |
|---|---|
| mAb | Monoclonal antibody |
| NSD | Nucleotide sugar donor |
| ER | Endoplasmic reticulum |
| QSS | Quasi-steady state |

| Abbreviations | Definition |
|---|---|
| DAE / PDAE / ODE | Differential algebraic equation / Partial differential algebraic equation / Ordinary differential equation |
| GDPMan, | Guanosine diphosphate mannose |
| GDPFuc, | Guanosine diphosphate fucose |
| UDPGal, | Uridine diphosphate galactose |
| UDPGlc, | Uridine diphosphate glucose |
| UDPGalNAc, | Uridine diphosphate N-acetylgalactosamine |
| UDPGlcNAc, | Uridine diphosphate N-acetylglucosamine |
| CMPNeu5Ac | Cytidine monophosphate N-acetylneuraminic acid |
| Amm, Asn, Asp, Glc, Gal, Gln, Glu, Lac, Urd | Ammonia, asparagine, aspartate, glucose, galactose, glutamine, glutamate, lactose, uridine) |
| Man5 | High-mannose glycan containing five mannose residues |
| FA2G0 | Core-fucosylated, biantennary complex glycan with no galactose residues |
| FA2G1 | Core-fucosylated, biantennary complex glycan with one galactose residue |
| FA2G2 | Core-fucosylated, biantennary complex glycan with two galactose residues—one on each antenna |
| G0 | Non-fucosylated, biantennary complex glycan with no galactose; shorthand for A2G0 |
| G2 | Non-fucosylated, biantennary complex glycan with two galactose units |
| OS / OS1 | Oligosaccharides / the first oligosaccharide |

| Variables | Description | Units |
|---|---|---|
| $q$ | Reaction rate in the cell culture | pg cell$^{-1}$ h$^{-1}$ for mAb and mmol cell$^{-1}$ h$^{-1}$ for the other metabolites |
| $f_{\mathrm{NSD}_i}$ | Flux of NSD consumed | µmol h$^{-1}$ |
| $Vel_{golgi}$ | Normalized linear velocity along the Golgi cisternae (PFR coordinate) | Golgi length min$^{-1}$ |
| $Y_i^{\mathrm{intra}}, Y_i^{\mathrm{extra}}$ | Intracellular / extracellular percentage of mAb bearing glycan $i$ | % |

| Subscripts | Description |
|---|---|
| $i$ | Metabolite / oligosaccharide index |
| $k$ | Nucleotide sugar donor index |
| $j$ | Reaction index |
| extra | Extracellular |
| intra | Intracellular |

| Superscripts | Description |
|---|---|
| glyc | Quantity specific to glycosylation reactions |
| hcp/lipid | Quantity specific to host cell protein and glycolipid synthesis |
| precursor | Quantity specific to precursor formation |